# QUENCHED INVARIANCE PRINCIPLE FOR MULTIDIMENSIONAL BALLISTIC RANDOM WALK IN A RANDOM ENVIRONMENT WITH A FORBIDDEN DIRECTION

By Firas Rassoul-Agha and Timo Seppäläinen[1]

*University of Utah and University of Wisconsin–Madison*

We consider a ballistic random walk in an i.i.d. random environment that does not allow retreating in a certain fixed direction. We prove an invariance principle (functional central limit theorem) under almost every fixed environment. The assumptions are nonnestling, at least two spatial dimensions, and a $2 + \varepsilon$ moment for the step of the walk uniformly in the environment. The main point behind the invariance principle is that the quenched mean of the walk behaves subdiffusively.

**1. Introduction.** This paper studies random walk in a random environment (RWRE) on the $d$-dimensional integer lattice $\mathbb{Z}^d$. This is a basic model in the field of disordered or random media. Our main result is a quenched invariance principle in dimension $d \geq 2$.

Here is a description of the model. An environment is a configuration of vectors of jump probabilities

$$\omega = (\omega_x)_{x \in \mathbb{Z}^d} \in \Omega = \mathcal{P}^{\mathbb{Z}^d},$$

where $\mathcal{P} = \{(p_z)_{z \in \mathbb{Z}^d} : p_z \geq 0, \sum_z p_z = 1\}$ is the simplex of all probability vectors on $\mathbb{Z}^d$. We use the notation $\omega_x = (\pi_{x,x+y})_{y \in \mathbb{Z}^d}$ for the coordinates of the probability vector $\omega_x$. The space $\Omega$ is equipped with the canonical product $\sigma$-field $\mathfrak{S}$ and with the natural shift $\pi_{xy}(T_z\omega) = \pi_{x+z,y+z}(\omega)$, for $z \in \mathbb{Z}^d$. On the space $(\Omega, \mathfrak{S})$ we are given an i.i.d. product measure $\mathbb{P}$. This means that the random probability vectors $(\omega_x)_{x \in \mathbb{Z}^d}$ are i.i.d. across the sites $x$ under $\mathbb{P}$.

Received September 2004; revised October 2005.

[1]Supported in part by NSF Grant DMS-04-02231.

The random walk operates as follows. An environment $\omega$ is chosen from the distribution $\mathbb{P}$ and fixed for all time. Pick an initial state $z \in \mathbb{Z}^d$. The random walk in environment $\omega$ started at $z$ is then the canonical Markov chain $\widehat{X} = (X_n)_{n \geq 0}$ with state space $\mathbb{Z}^d$ whose path measure $P_z^\omega$ satisfies

$$P_z^\omega(X_0 = z) = 1 \qquad \text{(initial state)},$$
$$P_z^\omega(X_{n+1} = y | X_n = x) = \pi_{xy}(\omega) \qquad \text{(transition probability)}.$$

The probability distribution $P_z^\omega$ on random walk paths is called the *quenched law*. The joint probability distribution

$$P_z(d\widehat{X}, d\omega) = P_z^\omega(d\widehat{X})\mathbb{P}(d\omega)$$

on walks and environments is called the *joint annealed law*, while its marginal on walks $P_z(d\widehat{X}, \Omega)$ is called simply the *annealed law*. $\mathbb{E}$, $E_0$ and $E_0^\omega$ denote expectations under, respectively, $\mathbb{P}$, $P_0$ and $P_0^\omega$.

We impose assumptions on the model that create a drift in some spatial direction $\hat{u}$. We also prohibit the walk from retreating in direction $\hat{u}$, a condition we express by saying that the walk has forbidden direction $-\hat{u}$. However, there is some freedom in the choice of $\hat{u}$. The long-term velocity $v$ of the walk need not be in direction $\hat{u}$, although of course the assumptions will imply $\hat{u} \cdot v > 0$.

We prove a quenched functional central limit theorem for the random walk. This means that, for $\mathbb{P}$-almost every $\omega$, under the measure $P_0^\omega$ the scaled random walk converges to a nondegenerate Brownian motion with a diffusion matrix that we describe. This result comes by a combination of regeneration, homogenization (studying the environment process) and martingale techniques. Our underlying proof strategy applies the approach of Maxwell and Woodroofe [9] and Derriennic and Lin [5] to the environment chain. This part is not spelled out in the present paper, but summarized in a theorem we quote from our earlier article [10]. The arguments of [9] and [5] themselves can be regarded as adaptations of the Kipnis–Varadhan method [8] to nonreversible situations.

The major technical part of our proof goes toward showing that the quenched mean $E_0^\omega(X_n)$ has variance of order $n^\gamma$ for some $\gamma < 1$. Bounding the variance of the quenched mean in turn is reduced to bounding the number of common points between two independent walks in a common environment. If we assume strictly more than a finite quenched third moment on the step of the walk, uniformly in the environment, we obtain $\gamma = 1/2$. Under a $p$th moment assumption with $2 < p \leq 3$ we can take any $\gamma > \frac{1}{p-1}$. The correct order of the variance of the quenched mean is an interesting question for this model, and for more general ballistic random walks. In the special case of space–time walks in $1+1$ dimensions with bounded steps, it has been proved that the quenched mean process, scaled by $n^{-1/4}$, converges to a certain Gaussian process [1].

The resulting quenched invariance principle admits two possible centerings, the asymptotic displacement $nv$ and the quenched mean. The approach and part of the result fail in one-dimensional walks and in certain other walks that are restricted to a single path by the environment (still considering only walks that satisfy the forbidden direction condition). In these cases a quenched invariance principle holds if the walk is centered at its quenched mean. But the quenched mean process itself also behaves diffusively with a Brownian motion limit. These other cases are explored in the paper [11].

There is a handful of quenched central limit theorems for RWRE in the literature. For the types of walks that we study, with a strong drift, Bolthausen and Sznitman [2] proved a quenched invariance principle. Their basic assumption is nonnestling which creates the drift, and for technical purposes they need an assumption of small noise and spatial dimension at least 4. (We get around these by making the forbidden direction assumption.) There is a certain analogy between our proof and the proof in [2]. Both proceed by bounding the variance of a certain quenched mean through control on the intersections of two independent paths. However, this similarity does not extend to the technical level, for we study a different variance and handle the intersections in a different manner.

For general overviews of recent developments in RWRE the reader can turn to the lectures [3, 12] and [13]. The introduction of [10] also presents a brief list of papers on central limit theorems for RWRE.

**2. Results.** Throughout the paper $\hat{u}$ is a fixed nonzero element of $\mathbb{R}^d$. We make a basic assumption called *nonnestling* that forces ballistic behavior on the walk.

HYPOTHESIS (N). There exists a positive deterministic constant $\delta$ such that
$$\mathbb{P}\left(\sum_z (z \cdot \hat{u}) \pi_{0z} \geq \delta\right) = 1.$$

In order to get the regeneration structure we need, we strengthen this assumption by requiring that the walk never retreats in the direction $\hat{u}$. Let us say the distribution $\mathbb{P}$ on environments has *forbidden direction* $-\hat{u}$ if

(2.1) $$\mathbb{P}\left(\sum_{z:z\cdot\hat{u}\geq 0} \pi_{0z} = 1\right) = 1.$$

This condition says that $X_n \cdot \hat{u}$ never decreases along the walk.

We also make a moment assumption uniformly in the environments. Let $|\cdot|$ denote the $\ell^1$- or the $\ell^2$-norm on $\mathbb{Z}^d$ (in all but a few computations the actual norm used is immaterial). For the invariance principle we need strictly more than a finite second moment, but other auxiliary results need fewer moments. Hence the power $p$ in the next hypothesis will be a parameter. Each time the hypothesis is invoked a bound for $p$ will be given, such as $p > 1$ or $p > 2$.

HYPOTHESIS (M).    There exist finite, deterministic positive constants $p$ and $M$ such that
$$\mathbb{P}\bigg(\sum_z |z|^p \pi_{0z} \leq M^p\bigg) = 1.$$

To take advantage of the renewal structure given by the nonnestling and forbidden direction assumptions, define a sequence of stopping times: $\sigma_0 = 0$, and for $k \geq 1$,

(2.2) $$\sigma_{k+1} = \inf\{n > \sigma_k : X_n \cdot \hat{u} \geq X_{\sigma_k} \cdot \hat{u} + 1\}.$$

Under the above assumptions the companion paper [11] shows these facts: $E_0(\sigma_k) < \infty$ for all $k$, $X_{\sigma_1}$ has $\bar{p}$th moment under $P_0$ for any $1 \leq \bar{p} < p$, and the walk has a long-term velocity $v = E_0(X_{\sigma_1})/E_0(\sigma_1)$ in the sense that $n^{-1}X_n \to v$ $P_0$-almost surely. See Theorem 3.3 and Lemma 3.4 in [11].

For the invariance principle we consider two centerings, the long-term displacement $nv$ and the quenched mean $E_0^\omega(X_n)$. So we define two scaled processes. For $t \in \mathbb{R}_+$ let

$$B_n(t) = \frac{X_{[nt]} - [nt]v}{\sqrt{n}} \quad \text{and} \quad \widetilde{B}_n(t) = \frac{X_{[nt]} - E_0^\omega(X_{[nt]})}{\sqrt{n}}.$$

Here $[x] = \max\{n \in \mathbb{Z} : n \leq x\}$ for $x \in \mathbb{R}$. Let $D_{\mathbb{R}^d}([0,\infty))$ denote the space of right-continuous $\mathbb{R}^d$-valued paths with left limits, endowed with the usual Skorohod topology (see the standard theory in [7]). For $\omega \in \Omega$ let $Q_n^\omega$, respectively $\widetilde{Q}_n^\omega$, denote the distribution of $B_n$, respectively $\widetilde{B}_n$, induced by $P_0^\omega$ on the Borel sets of $D_{\mathbb{R}^d}([0,\infty))$.

A quenched invariance principle cannot hold unless the walk is random under a fixed environment. This and more is contained in our final assumption of *ellipticity*.

HYPOTHESIS (E).    One has

(2.3) $$\mathbb{P}(\forall\, z \neq 0 : \pi_{0,0} + \pi_{0z} < 1) > 0.$$

Moreover, the walk is not supported by any one-dimensional subspace. More precisely, if $\mathcal{J} = \{y \in \mathbb{Z}^d : \mathbb{E}(\pi_{0y}) > 0\}$ is the set of all points that are accessible from 0 with one jump, then $\mathcal{J}$ is not contained in any subspace of the kind $\mathbb{R}u = \{su : s \in \mathbb{R}\}$ for any $u \in \mathbb{R}^d$. In particular, this rules out the case $d = 1$.

Let $\Gamma^t$ denote the transpose of a vector or matrix $\Gamma$. An element of $\mathbb{R}^d$ is regarded as a $d \times 1$ matrix, or column vector. For a symmetric, nonnegative definite $d \times d$ matrix $\Gamma$, a *Brownian motion with diffusion matrix* $\Gamma$ is the $\mathbb{R}^d$-valued process $\{W(t) : t \geq 0\}$ such that $W(0) = 0$, $W$ has continuous paths, independent increments, and for $s < t$ the $d$-vector $W(t) - W(s)$ has Gaussian distribution with mean zero and covariance matrix $(t - s)\Gamma$. The matrix $\Gamma$ is *degenerate* in direction $\xi \in \mathbb{R}^d$ if $\xi^t \Gamma \xi = 0$. Equivalently, $\xi \cdot W(t) = 0$ almost surely.

The diffusion matrix of our limiting process is defined by

(2.4) $$\mathfrak{D} = \frac{E_0[(X_{\sigma_1} - v\sigma_1)(X_{\sigma_1} - v\sigma_1)^t]}{E_0[\sigma_1]}.$$

One can check that this matrix $\mathfrak{D}$ is degenerate in direction $u$ if, and only if, $u$ is orthogonal to the vector space spanned by $\{x - y : \mathbb{E}(\pi_{0x})\mathbb{E}(\pi_{0y}) > 0\}$ (Theorem 4.1 in [11]). Degeneracy in directions that are orthogonal to all $x - y$, where $x$ and $y$ range over admissible jumps, cannot be avoided. This can be seen from the simple example of a homogeneous random walk that chooses with equal probability between two jumps $a$ and $b$. The diffusion matrix is then $\frac{1}{4}(a - b)(a - b)^t$.

We can now state the main theorem.

THEOREM 2.1. *Let $d \geq 2$ and consider an i.i.d. product probability measure $\mathbb{P}$ on environments with a forbidden direction $-\hat{u} \in \mathbb{Q}^d \setminus \{0\}$ as in (2.1). Assume nonnestling* (N) *in direction $\hat{u}$, the moment hypothesis* (M) *with $p > 2$, and ellipticity* (E). *Then as $n \to \infty$, for $\mathbb{P}$-almost every $\omega$ the distributions $Q_n^\omega$ and $\widetilde{Q}_n^\omega$ both converge weakly to the distribution of a Brownian motion with diffusion matrix $\mathfrak{D}$. Furthermore, the two centerings are asymptotically indistinguishable:*

$$\lim_{n \to \infty} \max_{0 \leq s \leq t} |\widetilde{B}_n(s) - B_n(s)| = \lim_{n \to \infty} n^{-1/2} \max_{k \leq [nt]} |E_0^\omega(X_k) - kv| = 0$$

*for $\mathbb{P}$-almost every $\omega$.*

Note that we assumed for Theorem 2.1 that the vector $\hat{u}$ has rational coordinates. Hypotheses (N) and (2.1) are not affected if $\hat{u}$ is multiplied by a constant. Hence later in the proof we can assume that $\hat{u}$ has integer coordinates.

In the special case where the step distribution of the walk is finitely supported, it turns out that if there is any nonzero vector $\hat{u}$ that satisfies both (2.1) and nonnestling (N), then there is also a rational one. We show this in Lemma A.1 in the Appendix. Thus for this case there is no restriction on $\hat{u}$. Since this case is perhaps the most important, we state it as a corollary.

COROLLARY 2.2. *Let $d \geq 2$ and consider an i.i.d. product probability measure $\mathbb{P}$ on environments with a forbidden direction $-\hat{u} \in \mathbb{R}^d \setminus \{0\}$ as in (2.1). Assume the step distribution is finitely supported, in other words that the set $\mathcal{J} = \{y \in \mathbb{Z}^d : \mathbb{E}(\pi_{0y}) > 0\}$ is finite. Assume nonnestling (N) in direction $\hat{u}$ and ellipticity (E). Then all the conclusions of Theorem 2.1 hold.*

We make a few remarks about ellipticity hypotheses. When (2.3) is violated, the environment $\omega$ determines completely the *set* of points $\{X_n : n \geq 0\}$ visited by the walk, and only the rate of advance remains random. In this case the process $B_n$ does not satisfy the quenched invariance principle. Same is true for the one-dimensional case. But $\widetilde{B}_n$ does satisfy an invariance principle in these cases, and furthermore, the quenched mean behaves diffusively. The companion paper [11] addresses these points.

One of the most popular hypotheses used in studies of RWRE is *uniform ellipticity*. One fixes a finite set $\mathcal{N}$ and a constant $0 < \kappa < 1$, and then assumes that $\mathbb{P}$-almost surely

$$\pi_{0z} = 0 \quad \text{for } z \neq \mathcal{N} \quad \text{and} \quad \kappa \leq \pi_{0z} \leq 1 \quad \text{for } z \in \mathcal{N}.$$

Suppose our forbidden direction assumption is made. Assume that $\mathcal{N}$ contains at least one point $x$ such that $x \cdot \hat{u} > 0$ and at least one other point $y$ such that $x$ and $y$ do not lie on a common line through the origin. Then all our other hypotheses (N), (M) and (E) follow. In particular, under the forbidden direction assumption, uniform ellipticity with a reasonably chosen $\mathcal{N}$ (such as that part of an $\ell^p$-ball of radius $\geq 1$ that satisfies $x \cdot \hat{u} \geq 0$) implies nonnestling.

The remainder of the paper proves Theorem 2.1. In several of our lemmas we indicate explicitly which assumptions are needed. In particular, $d \geq 2$ is not required everywhere, nor is the ellipticity assumption (E). We rely on a companion paper [11] for some basic results.

After the preliminaries the main work of the paper goes toward bounding the variance of the quenched mean. We record the result here.

THEOREM 2.3. *Let $d \geq 2$ and consider an i.i.d. product probability measure $\mathbb{P}$ on environments with a forbidden direction $-\hat{u} \in \mathbb{Q}^d \setminus \{0\}$ as in (2.1). Assume nonnestling (N) in direction $\hat{u}$, the moment hypothesis (M) with $p > 2$, and ellipticity (E). Let $\gamma > \frac{1}{p-1}$. Then there is a constant $C$ such that, for all $n \geq 1$,*

$$(2.5) \qquad \mathbb{E}[|E_0^\omega(X_n) - E_0(X_n)|^2] \leq \begin{cases} Cn^{1/2}, & \text{if } p > 3, \\ Cn^\gamma, & \text{if } 2 < p \leq 3. \end{cases}$$

Without affecting the validity of the bound (2.5), one can perform either one or both of these replacements: $E_0(X_n)$ can be replaced by $nv$, and $\mathbb{E}$ can be replaced by $\mathbb{E}_\infty$, expectation under the equilibrium measure of the environment chain introduced below in Theorem RS2. As pointed out in the Introduction, $n^{1/2}$ is known to be the correct order of the variance for some walks in $d = 2$.

**3. Preliminaries for the proof.** To prove the invariance principle we use the point of view of the particle. More precisely, we consider the Markov process on $\Omega$ with transition kernel

$$\hat{\pi}(\omega, A) = P_0^\omega(T_{X_1}\omega \in A).$$

For integers $n$ define $\sigma$-algebras $\mathfrak{S}_n = \sigma(\omega_x : x \cdot \hat{u} \geq n)$. Define the *drift* as

$$D(\omega) = E_0^\omega(X_1) = \sum_z z \pi_{0z}(\omega).$$

The proof of the quenched invariance principle is based on the next theorem from our earlier article [10]. This theorem is an application of the results of Maxwell and Woodroofe [9] and Derriennic and Lin [5] to random walk in random environment.

THEOREM RS1. *Let $d \geq 1$ and let $\mathbb{P}_\infty$ be any probability measure on $(\Omega, \mathfrak{S})$ that is invariant and ergodic for the Markov process on $\Omega$ with transition kernel $\hat{\pi}$. Assume that*

$$(3.1) \qquad \sum_z |z|^2 \mathbb{E}_\infty(\pi_{0z}) < \infty.$$

*Assume also that there exists an $0 \leq \alpha < 1/2$ such that*

$$(3.2) \qquad \mathbb{E}_\infty[|E_0^\omega(X_n) - n\mathbb{E}_\infty(D)|^2] = \mathcal{O}(n^{2\alpha}).$$

*Then for $\mathbb{P}_\infty$-almost every $\omega$ the distribution $Q_n^\omega$ of the process $\{B_n(t) : t \in \mathbb{R}_+\}$ converges weakly on the space $D_{\mathbb{R}^d}([0, \infty))$ to the distribution of a Brownian motion with a symmetric, nonnegative definite diffusion matrix that does not depend on $\omega$. Moreover, for $\mathbb{P}_\infty$-almost every $\omega$,*

$$(3.3) \qquad \lim_{n \to \infty} n^{-1/2} \max_{k \leq n} |E_0^\omega(X_k) - k\mathbb{E}_\infty(D)| = 0$$

*and, therefore, the same invariance principle holds for $\widetilde{Q}_n^\omega$.*

Above, $\mathbb{E}_\infty$ denotes expectation under the measure $\mathbb{P}_\infty$. To apply Theorem RS1, we need some preliminary results on equilibrium, the law of large numbers and the annealed invariance principle. These are contained in the next theorem that summarizes results from [11].

THEOREM RS2. *Let $d \geq 1$ and consider a product probability measure $\mathbb{P}$ on environments with a forbidden direction $-\hat{u} \in \mathbb{R}^d \setminus \{0\}$ as in (2.1). Assume nonnestling (N) in direction $\hat{u}$.*

(a) Ergodic equilibrium. *Assume the moment hypothesis* (M) *with $p > 1$. Then there exists a probability measure $\mathbb{P}_\infty$ on $(\Omega, \mathfrak{S})$ that is invariant for the Markov process with transition kernel $\hat{\pi}$ and has these properties:*

(i) $\mathbb{P} = \mathbb{P}_\infty$ *on $\mathfrak{S}_1$, $\mathbb{P}$ and $\mathbb{P}_\infty$ are mutually absolutely continuous on $\mathfrak{S}_0$, and $\mathbb{P}_\infty$ is absolutely continuous relative to $\mathbb{P}$ on $\mathfrak{S}_k$ with $k \leq 0$.*

(ii) *The Markov process with kernel $\hat{\pi}$ and initial distribution $\mathbb{P}_\infty$ is ergodic.*

(b) Law of large numbers. *Assume the moment hypothesis* (M) *with $p > 1$. Define $v = \mathbb{E}_\infty(D)$. Then we have the law of large numbers*

$$P_0\Big(\lim_{n\to\infty} n^{-1} X_n = v\Big) = 1.$$

*Moreover, $E_0(\sigma_1) < \infty$, $v = E_0(X_{\sigma_1})/E_0(\sigma_1)$ and*

(3.4) $$\sup_n |E_0(X_n) - nv| < \infty.$$

(c) Annealed invariance principle. *Assume the moment hypothesis* (M) *with $p > 2$. Then the distribution of the process $\{B_n(t) : t \in \mathbb{R}_+\}$ under $P_0$ converges weakly to the distribution of a Brownian motion with diffusion matrix $\mathfrak{D}$ defined by (2.4).*

The main idea for the proof of Theorem RS2 is that $(X_{\sigma_k} - X_{\sigma_{k-1}}, \sigma_k - \sigma_{k-1})_{k \geq 1}$ is a sequence of i.i.d. random variables under the annealed measure $P_0$.

Some comments follow. We have an explicit formula for the equilibrium distribution: if $A$ is $\mathfrak{S}_{-k}$-measurable for some $k \geq 0$, then

(3.5) $$\mathbb{P}_\infty(A) = \frac{E_0(\sum_{m=\sigma_k}^{\sigma_{k+1}-1} \mathbb{1}\{T_{X_m}\omega \in A\})}{E_0(\sigma_1)}.$$

The absolute continuity of $\mathbb{P}_\infty$ relative to $\mathbb{P}$ given by part (a) of Theorem RS2 has this consequence: moment assumption (M) is also valid under $\mathbb{P}_\infty$. Hence the drift $D$ can be integrated to define $v = \mathbb{E}_\infty(D)$. Then, if assumption (M) is strengthened to $p \geq 2$, it follows that hypothesis (3.1) of Theorem RS1 is fulfilled.

The course of the proof of Theorem 2.1 is now clear. Part (a) of Theorem RS2 gives the invariant measure needed for Theorem RS1. The real work goes toward checking (3.2). We first show, in Proposition 4.1 of Section 4, that it is enough to check (3.2) for $\mathbb{E}$ instead of $\mathbb{E}_\infty$. Then, in Sections 4 and 5, we check the latter condition is satisfied. At this point our proof will require more than two moments for $X_1$.

Suppose the hypotheses of Theorem RS1 have been checked. Let

$\mathcal{A} = \{\omega : Q_n^\omega \text{ and } \widetilde{Q}_n^\omega \text{ converge to the law of a Brownian motion and (3.3) holds}\}.$

The conclusion of Theorem RS1 is then $\mathbb{P}_\infty(\mathcal{A}) = 1$. Since $\mathcal{A}$ is $\mathfrak{S}_0$-measurable, mutual absolute continuity of $\mathbb{P}$ and $\mathbb{P}_\infty$ on $\mathfrak{S}_0$ implies that $\mathbb{P}(\mathcal{A}) = 1$. Theorem RS1 does not give the expression for the diffusion matrix. But the $\mathbb{P}$-almost sure quenched invariance principle must have the same limit as the annealed invariance principle. Hence part (c) of Theorem RS2 allows us to identify the limiting Brownian motion as the one with diffusion matrix $\mathfrak{D}$ from (2.4).

The upshot of this discussion is that in order to prove Theorem 2.1 only (3.2) remains to be verified. We finish this section of preliminaries by quoting part of Lemma 3.1 from [11]. Its proof uses standard ideas.

LEMMA 3.1. *Let $d \geq 1$ and consider a $T$-invariant probability measure $\mathbb{P}$ on environments with a forbidden direction $-\hat{u} \in \mathbb{R}^d \setminus \{0\}$ as in (2.1). Assume nonnestling (N) in direction $\hat{u}$, and the moment hypothesis (M) with $p > 1$. Then there exist strictly positive, finite constants $\bar{C}_m(M,\delta,p)$, $\hat{C}_{\bar{p}}(M,\delta,p)$ and $\lambda_0(M,\delta,p)$ such that for all $x \in \mathbb{Z}^d$, $\lambda \in [0,\lambda_0]$, $n,m \geq 0$ and $\mathbb{P}$-a.e. $\omega$,*

$$(3.6) \qquad E_x^\omega(|X_m - x|^{\bar{p}}) \leq M^{\bar{p}} m^{\bar{p}} \qquad \text{for } 1 \leq \bar{p} \leq p,$$

$$(3.7) \qquad P_x^\omega(\sigma_1 > n) \leq e^\lambda (1 - \lambda \delta/2)^n,$$

$$(3.8) \qquad E_x^\omega(\sigma_1^m) \leq \bar{C}_m,$$

$$(3.9) \qquad E_x^\omega(|X_{\sigma_1} - x|^{\bar{p}}) \leq \hat{C}_{\bar{p}} \qquad \text{for } 1 \leq \bar{p} < p.$$

**4. Bound for the variance of the quenched mean.** By the discussion in the previous section, it only remains to check (3.2) to derive the invariance principle Theorem 2.1 through an application of Theorem RS1. First we show in the next proposition that (3.2) is satisfied if it is true when $\mathbb{P}_\infty$ is replaced by $\mathbb{P}$. Subsequently we reduce this estimate to bounding the number of common points between two independent walks in a common environment. $X_{[0,n]} = \{X_k : 0 \leq k \leq n\}$ will denote the set of sites visited by the walk.

PROPOSITION 4.1. *Let $d \geq 1$ and consider a product probability measure $\mathbb{P}$ on environments with a forbidden direction $-\hat{u} \in \mathbb{R}^d \setminus \{0\}$ as in (2.1). Assume nonnestling (N) in direction $\hat{u}$, and the moment hypothesis (M) with $p \geq 2$. Let $\mathbb{P}_\infty$ be the measure in Theorem RS2(a). Assume that there exists an $\alpha < 1/2$ such that*

$$(4.1) \qquad \mathbb{E}(|E_0^\omega(X_n) - E_0(X_n)|^2) = \mathcal{O}(n^{2\alpha}).$$

*Then condition (3.2) is satisfied with this same $\alpha$.*

PROOF. Due to (3.4) the hypothesis becomes

$$\mathbb{E}(|E_0^\omega(X_n) - nv|^2) = \mathcal{O}(n^{2\alpha}). \tag{4.2}$$

Next, notice that $|v|^2 \leq M^2$ due to the moment hypothesis (M) and that $\mathbb{P}_\infty \ll \mathbb{P}$ on $\mathfrak{S}_0$. Notice also that (3.5) with $k=0$ implies that the Radon–Nikodym derivative $g_0 = d(\mathbb{P}_{\infty|\mathfrak{S}_0})/d(\mathbb{P}_{|\mathfrak{S}_0})$ is $\sigma(\omega_x : x \cdot \hat{u} < 1)$-measurable. Now the bound comes from a multistep calculation:

$$\mathbb{E}_\infty[|E_0^\omega(X_n) - n\mathbb{E}_\infty(D)|^2]$$
$$= \mathbb{E}_\infty[|E_0^\omega(X_n - nv, \sigma_1 \leq n) + E_0^\omega(X_n - nv, \sigma_1 > n)|^2]$$
$$\leq 2\mathbb{E}_\infty[|E_0^\omega\{X_n - (n-\sigma_1)v, \sigma_1 \leq n\} - E_0^\omega\{\sigma_1 v, \sigma_1 \leq n\}|^2]$$
$$+ 4(M^2 + |v|^2)n^2\mathbb{E}_\infty[P_0^\omega(\sigma_1 > n)^2]$$

by an application of (3.6),

$$\leq 4\mathbb{E}_\infty\left[\left|\sum_{x, m \leq n} P_0^\omega(X_m = x, \sigma_1 = m)E_x^\omega\{X_{n-m} - (n-m)v\}\right|^2\right]$$
$$+ 4(M^2 + 2|v|^2)\mathbb{E}_\infty[E_0^\omega(\sigma_1)^2]$$

by restarting the walk at time $\sigma_1$, by $|a+b|^2 \leq 2|a|^2 + 2|b|^2$ and by combining the expectations of $\sigma_1$,

$$\leq 4 \sum_{x, m \leq n} \mathbb{E}_\infty[P_0^\omega(X_m = x, \sigma_1 = m)|E_x^\omega\{X_{n-m} - (n-m)v\}|^2]$$
$$+ 12M^2 \mathbb{E}_\infty[E_0^\omega(\sigma_1)^2]$$

by an application of Jensen's inequality on the first term and by $|v| \leq M$,

$$= 4 \sum_{x, m \leq n} \mathbb{E}[g_0 P_0^\omega(X_m = x, \sigma_1 = m)|E_x^\omega\{X_{n-m} - (n-m)v\}|^2]$$
$$+ 12M^2 \mathbb{E}_\infty[E_0^\omega(\sigma_1)^2]$$
$$= 4 \sum_{x, m \leq n} \mathbb{E}[g_0 P_0^\omega(X_m = x, \sigma_1 = m)]\mathbb{E}[|E_x^\omega\{X_{n-m} - (n-m)v\}|^2]$$
$$+ 12M^2 \mathbb{E}_\infty[E_0^\omega(\sigma_1)^2]$$

because the i.i.d. assumption on $\mathbb{P}$ makes the two integrands independent,

$$\leq 8 \sum_{x, m \leq n} \mathbb{E}_\infty[P_0^\omega(X_m = x, \sigma_1 = m)]\mathbb{E}[|E_0^\omega\{X_{n-m} - (n-m)v\}|^2]$$
$$+ 8\mathbb{E}_\infty[E_0^\omega(|X_{\sigma_1}|^2)] + 12M^2 \mathbb{E}_\infty[E_0^\omega(\sigma_1)^2]$$

by shifting the initial state of $E_x^\omega$ back to 0, and by $|a+b|^2 \leq 2|a|^2 + 2|b|^2$ again,

$$= \mathcal{O}(n^{2\alpha}).$$

The final estimate above comes from (4.2) and the bounds in Lemma 3.1. □

Now, we will concentrate our attention on showing that (4.1) holds. This will be carried out in several steps. First, using Lemma 4.2, we bound the left-hand side of (4.1) by the expected number of intersections of two independent random walkers driven by the same environment. This is done in Proposition 4.3. Then, in Proposition 5.1 of Section 5, we bound this number of intersections and conclude the proof of Theorem 2.1.

For $U \subset \mathbb{Z}^d$ we use the notation $\omega_U = (\omega_x)_{x \in U}$. Recall that $X_{[0,n-1]}$ denotes the set of sites visited by the walk during time $0, \ldots, n-1$.

LEMMA 4.2. *Let $d \geq 1$ and consider a product probability measure $\mathbb{P}$ on environments with a forbidden direction $-\hat{u} \in \mathbb{R}^d \setminus \{0\}$ as in (2.1). Assume nonnestling* (N) *in direction $\hat{u}$, and the moment hypothesis* (M) *with $p > 1$. Fix $z \in \mathbb{Z}^d$ such that $z \cdot \hat{u} \geq 0$. Define the half-space $U = \{x \in \mathbb{Z}^d : x \cdot \hat{u} > z \cdot \hat{u}\}$. Let $\omega$ be an environment and $\tilde{\omega}$ another environment such that $\tilde{\omega}_x = \omega_x$ for all $x \neq z$. Then there exists a constant $C_0 = C_0(M, \delta)$ such that for all $z$, $\mathbb{P}$-almost every $\omega$, $\mathbb{P}$-almost every choice of $\tilde{\omega}_z$, and all $n \geq 1$,*

$$\left| \int [E_0^\omega(X_n) - E_0^{\tilde{\omega}}(X_n)] \mathbb{P}(d\omega_U) \right| \leq C_0 P_0^\omega(z \in X_{[0,n-1]}).$$

*Note that the right-hand side above is a function only of $\omega_{U^c}$ so there is no inconsistency.*

PROOF. Let $X_n$ and $\widetilde{X}_n$ denote walks that obey environments $\omega$ and $\tilde{\omega}$, respectively. We couple $X_n$ and $\widetilde{X}_n$ as follows. Given $\omega$, for each $x \in \mathbb{Z}^d$ pick a sequence of i.i.d. directed edges $(b_i(x) = (x, y_i))_{i \geq 1}$ from the distribution $(\pi_{xy}(\omega))_y$. Each time $X_n$ visits $x$, the walker takes a new edge $b_i(x)$, follows it to the next site $y_i$, discards the edge $b_i(x)$, and repeats this step at its new location. Since the edge $b_i(x)$ is discarded, next time $X_n$ visits $x$, $b_{i+1}(x)$ will be used.

The directed edges $\tilde{b}_i(x)$ that govern the walk $\widetilde{X}_n$ are defined by taking $\tilde{b}_i(x) = b_i(x)$ for $x \neq z$ and by picking i.i.d. directed edges $(\tilde{b}_i(z) = (z, y_i))_{i \geq 1}$ from the distribution $(\pi_{zy}(\tilde{\omega}))_y$.

Let $P_{x,\tilde{x}}^{\omega,\tilde{\omega}}$ denote this coupling measure under which the walks start at $x$ and $\tilde{x}$. If the walks start at 0, they stay together until they hit $z$. Let

$$\tau = \inf\{n \geq 0 : X_n = z\} = \inf\{n \geq 0 : \widetilde{X}_n = z\}$$

be the common hitting time of $z$ for the walks. Let

$$\sigma = \inf\{n \geq 0 : X_n \cdot \hat{u} > z \cdot \hat{u}\} \quad \text{and} \quad \tilde{\sigma} = \inf\{n \geq 0 : \widetilde{X}_n \cdot \hat{u} > z \cdot \hat{u}\}$$

be the times to enter the half-space $U$.

Note that $\sigma$ and $\tilde{\sigma}$ are different from $\sigma_1$ for a walk started at $z$. In fact, $\sigma \leq \sigma_1$.

We have

$$\begin{aligned} E_0^\omega(X_n) - E_0^{\tilde{\omega}}(X_n) &= E_{0,0}^{\omega,\tilde{\omega}}(X_n - \widetilde{X}_n) \\ &= E_{0,0}^{\omega,\tilde{\omega}}((X_n - \widetilde{X}_n)\mathbb{1}\{\tau < n\}) \\ &= E_{0,0}^{\omega,\tilde{\omega}}(X_n\mathbb{1}\{\tau < n\}) - E_{0,0}^{\omega,\tilde{\omega}}(\widetilde{X}_n\mathbb{1}\{\tau < n\}). \end{aligned}$$

Using the Markov property, one writes

$$E_{0,0}^{\omega,\tilde{\omega}}(X_n\mathbb{1}\{\tau < n\}) = \sum_{m=1}^n \sum_y P_{0,0}^{\omega,\tilde{\omega}}(\tau < n, \sigma \wedge n = m, X_m = y) E_y^\omega(X_{n-m}).$$

Note above that if $\tau < n$, then necessarily $\tau < \sigma \wedge n$, so the event $\{\tau < n, \sigma \wedge n = m\}$ is measurable with respect to $\sigma\{X_0, \ldots, X_m\}$. Rewrite the above as

$$\begin{aligned} E_{0,0}^{\omega,\tilde{\omega}}&(X_n\mathbb{1}\{\tau < n\}) \\ &= \sum_{1 \leq m, \tilde{m} \leq n} \sum_{y,\tilde{y}} P_{0,0}^{\omega,\tilde{\omega}}(\tau < n, \sigma \wedge n = m, \\ &\qquad\qquad \tilde{\sigma} \wedge n = \tilde{m}, X_m = y, \widetilde{X}_{\tilde{m}} = \tilde{y}) E_y^\omega(X_{n-m}). \end{aligned}$$

Develop the corresponding formula for $\widetilde{X}_n$, and subtract the two formulae to get

$$\begin{aligned} & E_0^\omega(X_n) - E_0^{\tilde{\omega}}(X_n) \\ (4.3) \quad &= \sum_{1 \leq m,\tilde{m} \leq n} \sum_{y,\tilde{y}} P_{0,0}^{\omega,\tilde{\omega}}(\tau < n, \sigma \wedge n = m, \tilde{\sigma} \wedge n = \tilde{m}, X_m = y, \widetilde{X}_{\tilde{m}} = \tilde{y}) \\ (4.4) \quad &\qquad \times (E_y^\omega(X_{n-m}) - E_{\tilde{y}}^\omega(\widetilde{X}_{n-\tilde{m}})). \end{aligned}$$

Note that the expectations on line (4.4) depend only on $\omega_U$. For $\widetilde{X}$ this is because if $\tilde{m} = n$, then $E_{\tilde{y}}^\omega(\widetilde{X}_{n-\tilde{m}}) = \tilde{y}$, while if $\tilde{m} < n$, then $\tilde{\sigma} = \tilde{m}$ and $\tilde{y} \in U$ and the walk never leaves $U$. The same reasoning works for the expectation of $X_{n-m}$. In fact, on line (4.4) we can drop the notational distinction between $X$ and $\widetilde{X}$.

Furthermore, the probabilities on line (4.3) are independent of $\omega_U$. They depend only on the environment in the complementary half-space $\{x \in \mathbb{Z}^d : x \cdot \hat{u} \leq z \cdot \hat{u}\}$.

Consider those terms in the sum on lines (4.3) and (4.4) with $m \leq \tilde{m}$. Then $n - m \geq n - \tilde{m}$ and we can write

$$E_y^\omega(X_{n-m}) - E_{\tilde{y}}^\omega(X_{n-\tilde{m}})$$
$$= E_y^\omega(X_{n-\tilde{m}}) - E_{\tilde{y}}^\omega(X_{n-\tilde{m}}) + E_y^\omega(X_{n-m} - X_{n-\tilde{m}})$$
$$= y - \tilde{y} + E_0^{T_y\omega}(X_{n-\tilde{m}}) - E_0^{T_{\tilde{y}}\omega}(X_{n-\tilde{m}}) + E_y^\omega(X_{n-m} - X_{n-\tilde{m}}).$$

Similarly for $m > \tilde{m}$,

$$E_y^\omega(X_{n-m}) - E_{\tilde{y}}^\omega(X_{n-\tilde{m}})$$
$$= y - \tilde{y} + E_0^{T_y\omega}(X_{n-m}) - E_0^{T_{\tilde{y}}\omega}(X_{n-m}) - E_{\tilde{y}}^\omega(X_{n-\tilde{m}} - X_{n-m}).$$

In both cases, when we integrate against $\mathbb{P}(d\omega_U)$ and use (3.6), we get

$$y - \tilde{y} + \text{(a term bounded in vector norm by } M|m - \tilde{m}|).$$

Substituting back into (4.3) and (4.4) gives

$$\left| \int (E_0^\omega(X_n) - E_0^{\tilde{\omega}}(X_n)) \mathbb{P}(d\omega_U) \right|$$
$$\leq \sum_{1 \leq m, \tilde{m} \leq n} \sum_{y, \tilde{y}} P_{0,0}^{\omega, \tilde{\omega}}(\tau < n, \sigma \wedge n = m, \tilde{\sigma} \wedge n = \tilde{m},$$
$$X_m = y, \widetilde{X}_{\tilde{m}} = \tilde{y}) |y - \tilde{y}|$$
$$+ \sum_{1 \leq m, \tilde{m} \leq n} \sum_{y, \tilde{y}} P_{0,0}^{\omega, \tilde{\omega}}(\tau < n, \sigma \wedge n = m, \tilde{\sigma} \wedge n = \tilde{m},$$
$$X_m = y, \widetilde{X}_{\tilde{m}} = \tilde{y}) M|m - \tilde{m}|.$$

For the first term to the right-hand side of the inequality, noting again that on the event $\{\tau < n\}$ we have $\tau < \sigma \wedge \tilde{\sigma} \wedge n$, one can write

$$\sum_{1 \leq m, \tilde{m} \leq n} \sum_{y, \tilde{y}} P_{0,0}^{\omega, \tilde{\omega}}(\tau < n, \sigma \wedge n = m, \tilde{\sigma} \wedge n = \tilde{m}, X_m = y, \widetilde{X}_{\tilde{m}} = \tilde{y}) |y - \tilde{y}|$$
$$= E_{0,0}^{\omega, \tilde{\omega}}(\mathbb{1}\{\tau < n\} |X_{\sigma \wedge n} - \widetilde{X}_{\tilde{\sigma} \wedge n}|)$$
(4.5) $$= \sum_{\ell=0}^{n-1} E_{0,0}^{\omega, \tilde{\omega}}(\mathbb{1}\{\tau = \ell\} |X_{\sigma \wedge n} - \widetilde{X}_{\tilde{\sigma} \wedge n}|)$$
$$= \sum_{\ell=0}^{n-1} P_0^\omega(\tau = \ell) E_{z,z}^{\omega, \tilde{\omega}}(|X_{\sigma \wedge (n-\ell)} - \widetilde{X}_{\tilde{\sigma} \wedge (n-\ell)}|)$$
$$\leq \sum_{\ell=0}^{n-1} P_0^\omega(\tau = \ell) [E_z^\omega(|X_{\sigma \wedge (n-\ell)} - z|) + E_z^{\tilde{\omega}}(|X_{\sigma \wedge (n-\ell)} - z|)].$$

Note now that by (3.6) and (3.7) we have, for all $n \geq 0$,

$$\begin{aligned}
E_z^\omega(|X_{\sigma \wedge n} - z|) &\leq E_z^\omega(|X_\sigma - z|) + E_z^\omega(|X_n - z|, \sigma > n) \\
&\leq \sum_{m \geq 1} E_z^\omega(|X_m - z|^p)^{1/p} P_z^\omega(X_{m-1} \cdot \hat{u} = z \cdot \hat{u})^{(p-1)/p} \\
&\quad + E_z^\omega(|X_n - z|^p)^{1/p} P_z^\omega(X_n \cdot \hat{u} = z \cdot \hat{u})^{(p-1)/p} \\
&\leq \hat{C}.
\end{aligned}$$

Therefore, we can bound (4.5) by $\hat{C} P_0^\omega(\tau < n)$. For the remaining sum there is a similar bound:

$$\begin{aligned}
\sum_{1 \leq m, \tilde{m} \leq n} \sum_{y, \tilde{y}} &P_{0,0}^{\omega, \tilde{\omega}}(\tau < n, \sigma \wedge n = m, \tilde{\sigma} \wedge n = \tilde{m}, X_m = y, \widetilde{X}_{\tilde{m}} = \tilde{y}) M|m - \tilde{m}| \\
&= M E_{0,0}^{\omega, \tilde{\omega}}(\mathbb{1}\{\tau < n\}|\sigma \wedge n - \tilde{\sigma} \wedge n|) \\
&= M \sum_{\ell = 0}^{n-1} P_0^\omega(\tau = \ell) E_{z,z}^{\omega, \tilde{\omega}}[|\sigma \wedge (n-\ell) - \tilde{\sigma} \wedge (n-\ell)|] \\
&\leq M P_0^\omega(\tau < n)[E_z^\omega(\sigma_1) + E_z^{\tilde{\omega}}(\sigma_1)] \\
&\leq 2\bar{C}_1 M P_0^\omega(\tau < n).
\end{aligned}$$

Putting the bounds together gives

$$\left| \int [E_0^\omega(X_n) - E_0^{\tilde{\omega}}(X_n)] \mathbb{P}(d\omega_U) \right| \leq (\hat{C} + 2\bar{C}_1 M) P_0^\omega(\tau < n),$$

which is the claim. □

Now we take one step toward proving (4.1). We write $P_{x,y}$ and $E_{x,y}$ for probabilities and expectations on a probability space on which are defined the $\mathbb{P}$-distributed environments, and two walks $X_n$ and $\widetilde{X}_n$ that are independent given the environment, and whose initial points are $X_0 = x$ and $\widetilde{X}_0 = y$. Similarly, $P_{x,y}^\omega$ and $E_{x,y}^\omega$ will be the quenched probabilities and expectations. Note that this coupling of walks $X_n$ and $\widetilde{X}_n$ is quite different from the one in the proof of Lemma 4.2. Let $|A|$ denote the cardinality of a set $A \subseteq \mathbb{Z}^d$. We have the following:

PROPOSITION 4.3. *Let $d \geq 1$ and consider a product probability measure $\mathbb{P}$ on environments with a forbidden direction $-\hat{u} \in \mathbb{R}^d \setminus \{0\}$ as in (2.1). Assume nonnestling* (N) *in direction $\hat{u}$, and the moment hypothesis* (M) *with $p > 1$. Let $C_0$ be as in Lemma* 4.2. *Then we have for all $n \geq 0$,*

$$(4.6) \qquad \mathbb{E}[|E_0^\omega(X_n) - E_0(X_n)|^2] \leq C_0^2 E_{0,0}(|X_{[0,n-1]} \cap \widetilde{X}_{[0,n-1]}|).$$

PROOF. For $L \geq 0$, define $B_L = \{x \in \mathbb{Z}^d : |x| \leq L\}$. Also, for $B \subset \mathbb{Z}^d$, let $\mathfrak{S}_B = \sigma(\omega_B)$. Fix $n \geq 1$ and $L \geq 0$ and let $(x_j)_{j \geq 1}$ be some fixed ordering of $B_L$ satisfying

$$\forall i \geq j : x_i \cdot \hat{u} \geq x_j \cdot \hat{u}.$$

Set $\mathcal{U}_0$ to be the trivial $\sigma$-field and define the filtration $\mathcal{U}_j = \sigma(\omega_{x_1}, \ldots, \omega_{x_j})$ and the variables $\zeta_j = \mathbb{E}(E_0^\omega(X_n)|\mathcal{U}_j)$.

$(\zeta_j - \zeta_{j-1})_{j \geq 1}$ is a sequence of $L^2(\mathbb{P})$-martingale differences, and so

$$\mathbb{E}[|\mathbb{E}(E_0^\omega\{X_n\}|\mathfrak{S}_{B_L}) - E_0(X_n)|^2] = \sum_{j=1}^{|B_L|} \mathbb{E}(|\zeta_j - \zeta_{j-1}|^2)$$

$$\leq C_0^2 \sum_{z \in B_L} \mathbb{E}[P_0^\omega(z \in X_{[0,n-1]})^2]$$

$$\leq C_0^2 \sum_z \mathbb{E}[P_{0,0}^\omega(z \in X_{[0,n-1]} \cap \widetilde{X}_{[0,n-1]})]$$

$$= C_0^2 \mathbb{E}[E_{0,0}^\omega(|X_{[0,n-1]} \cap \widetilde{X}_{[0,n-1]}|)],$$

where the first inequality is due to Lemma 4.2. By (3.6) $E_0^\omega(X_n)$ is a bounded random variable and, therefore, $\mathbb{E}[E_0^\omega(X_n)|\mathfrak{S}_{B_L}]$ converges in $L^2(\mathbb{P})$ to $E_0^\omega(X_n)$. Thus, taking $L$ to infinity proves the proposition. $\square$

**5. Bound for number of common points between two independent paths.** In this section we show that the right-hand side of (4.6) is $\mathcal{O}(n^{1-\delta})$ where $\delta > 0$ depends on the strength of our moment hypothesis (M). We say that $x$ belongs to *level* $\ell$ if $x \cdot \hat{u} = \ell$. We will count the number of common points between two paths by levels. This is where the assumption that $\hat{u}$ is a rational vector is needed. Otherwise the levels could accumulate and we would not be able to number them. As mentioned in the remarks following Theorem 2.1, if $\hat{u} \in \mathbb{Q}^d \setminus \{0\}$, then we can and will assume, without any loss of generality, that $\hat{u} \in \mathbb{Z}^d \setminus \{0\}$. This way we only need to consider integral levels $\ell$. The assumption of integral $\hat{u}$ also has the effect that the stopping times $\{\sigma_k\}$ defined by (2.2) mark the successive jumps to new levels. Define

$$\mathbb{V}_d = \{y \in \mathbb{Z}^d : y \cdot \hat{u} = 0\}$$

and recall that $\mathcal{J} = \{y : \mathbb{E}(\pi_{0y}) > 0\}$. Recall also that under $P_{x,y}^\omega$ the walks $X$ and $\widetilde{X}$ are independent in the common environment $\omega$ with initial points $X_0 = x$ and $\widetilde{X}_0 = y$, and $P_{x,y} = \int P_{x,y}^\omega \mathbb{P}(d\omega)$. Now for the first time we need the ellipticity assumptions.

PROPOSITION 5.1. *Let $d \geq 2$ and consider a product probability measure $\mathbb{P}$ on environments with a forbidden direction $-\hat{u} \in \mathbb{Z}^d \setminus \{0\}$ as in (2.1).*

*Assume nonnestling* (N) *in direction* $\hat{u}$, *the moment hypothesis* (M) *with* $p > 2$ *and ellipticity* (E). *Let* $\gamma > \frac{1}{p-1}$. *Then there exists a constant* $C_1 < \infty$ *such that*

$$E_{0,0}(|X_{[0,n-1]} \cap \widetilde{X}_{[0,n-1]}|) \leq \begin{cases} C_1 n^{1/2}, & \text{if } p > 3, \\ C_1 n^{\gamma}, & \text{if } 2 < p \leq 3. \end{cases}$$

PROOF. Denote the times of reaching a level at or above $\ell$ by

$$\lambda_\ell = \inf\{n \geq 0 : X_n \cdot \hat{u} \geq \ell\} \quad \text{and} \quad \tilde{\lambda}_\ell = \inf\{n \geq 0 : \widetilde{X}_n \cdot \hat{u} \geq \ell\}.$$

We may occasionally write $\lambda(\ell)$ for $\lambda_\ell$ to avoid complicated subscripts on subscripts. Note that $X$ hits level $\ell$ if, and only if, $X_{\lambda_\ell} \cdot \hat{u} = \ell$. Common points of the two paths can occur only on levels that are visited by both paths, or "common levels." These common levels are denoted by the random variables $0 = L_0 < L_1 < L_2 < \cdots$ defined by

$$L_j = \inf\{\ell > L_{j-1} : X_{\lambda_\ell} \cdot \hat{u} = \widetilde{X}_{\tilde{\lambda}_\ell} \cdot \hat{u} = \ell\}.$$

Let $\mathcal{F}_n$ be the filtration of the walk $X_n$, and similarly $\widetilde{\mathcal{F}}_n$ for $\widetilde{X}_n$. Let $\mathcal{H}_0$ be the trivial $\sigma$-field, and

$$\mathcal{H}_\ell = \sigma(\{\omega_x : x \cdot \hat{u} < \ell\}, \mathcal{F}_{\lambda_\ell}, \widetilde{\mathcal{F}}_{\tilde{\lambda}_\ell}).$$

The $L_j$'s are stopping times for the filtration $\{\mathcal{H}_\ell\}$. Lemma 5.3 below shows that $L_j$ is finite for all $j$.

Now we can rewrite the mean number of common points as follows. Write temporarily

$$N_\ell = |\{x \in \mathbb{Z}^d : x \cdot \hat{u} = \ell, x \in X_{[0,\infty)} \cap \widetilde{X}_{[0,\infty)}\}|$$

for the number of common points on level $\ell$:

$$E_{0,0}(|X_{[0,n-1]} \cap \widetilde{X}_{[0,n-1]}|)$$

$$\leq \sum_{\ell=0}^{\infty} E_{0,0}[N_\ell \mathbb{1}\{X_{\lambda_\ell} \cdot \hat{u} = \widetilde{X}_{\tilde{\lambda}_\ell} \cdot \hat{u} = \ell, \lambda_\ell \vee \tilde{\lambda}_\ell < n\}]$$

$$= \sum_{\ell=0}^{\infty} E_{0,0}[E_{0,0}(N_\ell | \mathcal{H}_\ell) \mathbb{1}\{X_{\lambda_\ell} \cdot \hat{u} = \widetilde{X}_{\tilde{\lambda}_\ell} \cdot \hat{u} = \ell, \lambda_\ell \vee \tilde{\lambda}_\ell < n\}].$$

Introduce the function

$$h(z) = E_{z,0}(|X_{[0,\sigma_1)} \cap \widetilde{X}_{[0,\tilde{\sigma}_1)}|)$$

for $z \in \mathbb{V}_d$. Then on the event $X_{\lambda_\ell} \cdot \hat{u} = \widetilde{X}_{\tilde{\lambda}_\ell} \cdot \hat{u} = \ell$

$$E_{0,0}(N_\ell | \mathcal{H}_\ell) = h(X_{\lambda_\ell} - \widetilde{X}_{\tilde{\lambda}_\ell}).$$

Introduce the process

(5.1) $$Z_j = X_{\lambda(L_j)} - \widetilde{X}_{\tilde{\lambda}(L_j)} \in \mathbb{V}_d$$

to rewrite the previous development as

$$E_{0,0}(|X_{[0,n-1]} \cap \widetilde{X}_{[0,n-1]}|)$$
$$\leq \sum_{\ell=0}^{\infty} E_{0,0}[h(X_{\lambda_\ell} - \widetilde{X}_{\tilde{\lambda}_\ell})\mathbb{1}\{X_{\lambda_\ell} \cdot \hat{u} = \widetilde{X}_{\tilde{\lambda}_\ell} \cdot \hat{u} = \ell, \lambda_\ell \vee \tilde{\lambda}_\ell < n\}]$$
$$= \sum_{j=0}^{\infty} E_{0,0}[h(Z_j)\mathbb{1}\{\lambda_{L_j} \vee \tilde{\lambda}_{L_j} < n\}].$$

Finally observe that $\lambda_{L_j} \geq j$ because it takes at least $j$ jumps to get to the $j$th common level. This gives us the inequality

(5.2) $$E_{0,0}(|X_{[0,n-1]} \cap \widetilde{X}_{[0,n-1]}|) \leq \sum_{j=0}^{n-1} E_0[h(Z_j)].$$

Equation (5.2) is the starting point for the analysis. The subscript in the last $E_0$ above is the initial point $Z_0 = 0 \in \mathbb{V}_d$.

To complete the proof of Proposition 5.1 we need to control the function $h$ and the process $Z_j$. We start with $h$.

LEMMA 5.2. *Let $d \geq 1$ and consider a product probability measure $\mathbb{P}$ on environments with a forbidden direction $-\hat{u} \in \mathbb{Z}^d \setminus \{0\}$ as in (2.1). Assume nonnestling (N) in direction $\hat{u}$, and the moment hypothesis (M) with $p > 2$. Then the function $h$ is summable:*

$$\sum_{z \in \mathbb{V}_d} h(z) < \infty.$$

PROOF. Define $b(x) = |x| + 1$ for $x \in \mathbb{Z}^d$. Below we will use the properties $b(x) = b(-x)$ and $b(x+y) \leq b(x)b(y)$. Notice that the number of points on the path $X_{[0,\sigma_1)}$ is at most $\sigma_1$. Bound $h(0)$ simply by $h(0) \leq E_0(\sigma_1)$. We bound the sum of the remaining terms as follows:

$$\sum_{z \cdot \hat{u}=0, z \neq 0} h(z) \leq \sum_{z \cdot \hat{u}=0, z \neq 0} E_{z,0}(\sigma_1 \mathbb{1}\{X_{[0,\sigma_1)} \cap \widetilde{X}_{[0,\tilde{\sigma}_1)} \neq \varnothing\})$$
$$\leq \sum_{z \cdot \hat{u}=0, z \neq 0} E_{z,0}(\sigma_1 \mathbb{1}\{\sigma_1 > b(z)\})$$
$$+ \sum_{z \cdot \hat{u}=0, z \neq 0} b(z) P_{z,0}(X_{[0,\sigma_1)} \cap \widetilde{X}_{[0,\tilde{\sigma}_1)} \neq \varnothing).$$

The first sum after the last inequality is finite by the exponential tail bounds (3.7).

We decompose the last sum according to the first site $y$ along the $X$-walk that is also visited by the $\widetilde{X}$-walk. Then the $\widetilde{X}$-walk from 0 to $y$ does not intersect the $X$-walk from $z$ to $y$, except at $y$. To formalize this, for $z \neq 0$ and $y \in \mathbb{V}_d$, let $\Gamma(z,0,y)$ be the set of all pairs of paths $(\gamma, \tilde{\gamma})$ such that $\gamma = \{z = x_0, x_1, \ldots, x_m = y\}$, $\tilde{\gamma} = \{0 = y_0, y_1, \ldots, y_n = y\}$, all points reside on level 0, and $y$ is the first common point along the two paths. Paths $\gamma = \{z\}$ and $\tilde{\gamma} = \{0\}$ are also considered when either $y = z$ or $y = 0$. Use the notation

$$P^\omega(\gamma) = \pi_{x_0,x_1}(\omega) \pi_{x_1,x_2}(\omega) \cdots \pi_{x_{m-1},x_m}(\omega)$$

for the probability that the $X$-walk follows path $\gamma$, and similarly for $P^\omega(\tilde{\gamma})$. For any pair $(\gamma, \tilde{\gamma}) \in \Gamma(z,0,y)$ the random variables $P^\omega(\gamma)$ and $P^\omega(\tilde{\gamma})$ are independent under $\mathbb{P}$. Let $H(z,y)$ be the collection of all paths from $z$ to $y$ on level 0 that contain $y$ only as the last site, and analogously for $H(0,y)$. Then

$$\sum_{z \cdot \hat{u} = 0, z \neq 0} b(z) P_{z,0}(X_{[0,\sigma_1)} \cap \widetilde{X}_{[0,\tilde{\sigma}_1)} \neq \varnothing)$$

$$= \sum_{z \cdot \hat{u} = 0, z \neq 0} b(z) \sum_{y \cdot \hat{u} = 0} \sum_{(\gamma,\tilde{\gamma}) \in \Gamma(z,0,y)} \mathbb{E}[P^\omega(\gamma)] \mathbb{E}[P^\omega(\tilde{\gamma})]$$

$$\leq \sum_{y \cdot \hat{u} = 0} b(y) \sum_{\tilde{\gamma} \in H(0,y)} \mathbb{E}[P^\omega(\tilde{\gamma})] \sum_{z \cdot \hat{u} = 0} b(y-z) \sum_{\gamma \in H(z,y)} \mathbb{E}[P^\omega(\gamma)]$$

$$= \sum_{y \cdot \hat{u} = 0} b(y) \sum_{\tilde{\gamma} \in H(0,y)} \mathbb{E}[P^\omega(\tilde{\gamma})] \sum_{x \cdot \hat{u} = 0} b(x) \sum_{\gamma \in H(0,x)} \mathbb{E}[P^\omega(\gamma)]$$

$$= \left( \sum_{x \cdot \hat{u} = 0} b(x) \sum_{\gamma \in H(0,x)} \mathbb{E}[P^\omega(\gamma)] \right)^2 = \left( \sum_{x \cdot \hat{u} = 0} b(x) P_0(x \in X_{[0,\sigma_1)}) \right)^2$$

$$\leq \left( E_0 \left[ \sum_{n=0}^{\sigma_1 - 1} (1 + |X_n|) \right] \right)^2 < \infty.$$

The finiteness of the last term can be seen by the usual application of Hölder's inequality to $E_0(|X_n| \mathbb{1}\{\sigma_1 > n\})$, along with (3.7) and (3.9).  $\square$

Next we analyze the process $Z_j$. Under the annealed probability it is a Markov chain because the walks can be restarted from the points $(X_{\lambda(L_j)}, \widetilde{X}_{\tilde{\lambda}(L_j)})$ of each new common level, and then they see a new environment independent of the past. We shall show that $Z_j$ is also a martingale with certain uniform moment bounds on its increments.

Let $L = L_1$ denote the first common level above zero. We generalize the treatment to two walks $X_n$ and $\widetilde{X}_n$ that both start at level zero, but not necessarily at the same point. The first task is to bound $L$.

LEMMA 5.3. *Let $d \geq 1$ and consider a product probability measure $\mathbb{P}$ on environments with a forbidden direction $-\hat{u} \in \mathbb{Z}^d \setminus \{0\}$ as in (2.1). Assume nonnestling (N) in direction $\hat{u}$, and the moment hypothesis (M) with $p > 2$. Then for any $\bar{p} \in [2, p)$ there exists a constant $C_2 = C_2(\bar{p})$ such that $E_{z,0}(L^{\bar{p}-1}) \leq C_2$ for all choices of $z \in \mathbb{V}_d$.*

PROOF. Either the very first new levels of $X$ and $\widetilde{X}$ are common, or not, so

$$
\begin{aligned}
E_{z,0}(L^{\bar{p}-1}) &= E_{z,0}[(X_{\sigma_1} \cdot \hat{u})^{\bar{p}-1} \mathbb{1}\{X_{\sigma_1} \cdot \hat{u} = \widetilde{X}_{\tilde{\sigma}_1} \cdot \hat{u}\}] \\
&\quad + \sum_{x \cdot \hat{u} \neq \tilde{x} \cdot \hat{u}} E_{z,0}[L^{\bar{p}-1} \mathbb{1}\{X_{\sigma_1} = x\} \mathbb{1}\{\widetilde{X}_{\tilde{\sigma}_1} = \tilde{x}\}].
\end{aligned}
\tag{5.3}
$$

The first term after the equality sign is bounded by a constant independently of $z$ by (3.9). We rewrite the last sum by introducing the levels visited by the $X$-walk until the first common level. It becomes

$$
\sum_{\substack{k \geq 1 \\ (\tilde{i}, i_1, \ldots, i_k) \in A_k}} \sum_{x \cdot \hat{u} = i_1, \tilde{x} \cdot \hat{u} = \tilde{i}} i_k^{\bar{p}-1} P_{z,0}(X_{\sigma_1} = x, X_{\sigma_m} \cdot \hat{u} = i_m \text{ for } m = 2, \ldots, k,
$$

$$
\widetilde{X}_{\tilde{\sigma}_1} = \tilde{x}, \widetilde{X}\text{-walk does not visit}
$$

$$
\text{levels } i_1, \ldots, i_{k-1} \text{ but does visit level } i_k),
$$

where $A_k$ is the set of positive integer $(k+1)$-vectors $(\tilde{i}, i_1, \ldots, i_k)$ such that:

  (i) if $k = 1$, then $0 < \tilde{i} < i_1$, while
  (ii) if $k \geq 2$, then $0 < i_1 < \cdots < i_k$, $\tilde{i} \leq i_k$ and $\tilde{i} \notin \{i_1, \ldots, i_{k-1}\}$.

This accounts for all the possible ways of saying that the walks continue from disjoint levels $i_1$ and $\tilde{i}$ and first meet at level $i_k$. It can happen that $i_k = i_1$ or $i_k = \tilde{i}$, but not both.

Write the probability in the above sum as

$$
\mathbb{E}[P_z^\omega(X_{\sigma_1} = x) P_0^\omega(\widetilde{X}_{\tilde{\sigma}_1} = \tilde{x}) P_x^\omega(X_{\sigma_m} \cdot \hat{u} = i_{m+1} \text{ for } m = 1, \ldots, k-1)
$$

$$
\times P_{\tilde{x}}^\omega(\widetilde{X}\text{-walk does not visit levels } i_1, \ldots, i_{k-1} \text{ but does visit level } i_k)]
$$

$$
= P_{z,0}(X_{\sigma_1} = x, \widetilde{X}_{\tilde{\sigma}_1} = \tilde{x}) P_x(X_{\sigma_m} \cdot \hat{u} = i_{m+1} \text{ for } m = 1, \ldots, k-1)
$$

$$
\times P_{\tilde{x}}(\widetilde{X}\text{-walk does not visit levels } i_1, \ldots, i_{k-1} \text{ but does visit level } i_k).
$$

Above we used independence: the probabilities

$$
P_z^\omega(X_{\sigma_1} = x) P_0^\omega(\widetilde{X}_{\tilde{\sigma}_1} = \tilde{x})
$$

are functions of $(\omega_y : y \cdot \hat{u} = 0)$, the probability

$$P_x^\omega(X_{\sigma_m} \cdot \hat{u} = i_{m+1} \text{ for } m = 1, \ldots, k-1)$$

is a function of $(\omega_y : y \cdot \hat{u} \in \{i_1, \ldots, i_{k-1}\})$, while probability

$$P_{\tilde{x}}^\omega(\widetilde{X}\text{-walk does not visit levels } i_1, \ldots, i_{k-1} \text{ but does visit level } i_k)$$

is a function of $(\omega_y : 0 < y \cdot \hat{u} < i_k, y \cdot \hat{u} \notin \{i_1, \ldots, i_{k-1}\})$.

By translation, the last sum in (5.3) can now be written as

$$\sum_{i_1 \neq \tilde{i}} P_{z,0}(X_{\sigma_1} \cdot \hat{u} = i_1, \widetilde{X}_{\tilde{\sigma}_1} \cdot \hat{u} = \tilde{i})$$

$$\times \left\{ i_1^{\bar{p}-1} P_0(\tilde{i} + \widetilde{X}_{\sigma_n} \cdot \hat{u} = i_1 \text{ for some } n \geq 1) \right.$$

$$+ \sum_{\substack{k \geq 2, (i_2, \ldots, i_k): \\ (\tilde{i}, i_1, \ldots, i_k) \in A_k}} i_k^{\bar{p}-1} P_0(i_1 + X_{\sigma_j} \cdot \hat{u} = i_{j+1} \text{ for } j = 1, \ldots, k-1)$$

$$\times P_0(\tilde{i} + \widetilde{X}_{\sigma_j} \cdot \hat{u} \notin \{i_1, \ldots, i_{k-1}\} \text{ for all } j \geq 0,$$

$$\left. \text{but } \tilde{i} + \widetilde{X}_{\sigma_n} \cdot \hat{u} = i_k \text{ for some } n \geq 0) \right\}.$$

The quantity in braces can be represented as $E(L_{i_1,\tilde{i}}^{\bar{p}-1})$ where the random variable $L_{i_1,\tilde{i}}$ is defined as the first common time (or "level") of two independent delayed renewal processes:

$$L_{i,j} = \inf\left\{ \ell \geq 1: \text{ for some } m, n \geq 0, i + \sum_{k=1}^m Y_k = \ell = j + \sum_{k=1}^n \widetilde{Y}_k \right\},$$

where $\{Y_k\}$ is an i.i.d. positive integer-valued sequence distributed like $\{(X_{\sigma_k} - X_{\sigma_{k-1}}) \cdot \hat{u}\}$ under $P_0$, and $\{\widetilde{Y}_k\}$ is an independent copy. It follows from Lemma 3.1 that $E(Y_1^{\bar{p}}) < \infty$. By Lemma A.3 in the Appendix, $E(L_{i_1,\tilde{i}}^{\bar{p}-1}) \leq C(1 + i_1^{\bar{p}-1} + \tilde{i}^{\bar{p}-1})$. Substituting this back into (5.3) gives

$$E_{z,0}(L^{\bar{p}-1}) \leq C + C \sum_{i_1 \neq \tilde{i}} P_{z,0}(X_{\sigma_1} \cdot \hat{u} = i_1, \widetilde{X}_{\tilde{\sigma}_1} \cdot \hat{u} = \tilde{i})(1 + i_1^{\bar{p}-1} + \tilde{i}^{\bar{p}-1})$$

which is bounded by a constant by (3.9). This completes the proof of Lemma 5.3.
□

Having bounded $L$, we turn to develop a martingale. We have

$$E_{z,0}^\omega(X_{\lambda_{k+1}}|\mathcal{H}_k) = \mathbb{1}\{X_{\lambda_k} \cdot \hat{u} \geq k+1\} X_{\lambda_k}$$

$$+ \mathbb{1}\{X_{\lambda_k} \cdot \hat{u} = k\}(X_{\lambda_k} + E^{\omega}_{X_{\lambda_k}}(X_{\sigma_1} - X_0))$$

$$= X_{\lambda_k} + \mathbb{1}\{X_{\lambda_k} \cdot \hat{u} = k\} E^{\omega}_{X_{\lambda_k}}(X_{\sigma_1} - X_0).$$

Consequently

$$M_k = X_{\lambda_k} - \sum_{j=0}^{k-1} \mathbb{1}\{X_{\lambda_j} \cdot \hat{u} = j\} E^{\omega}_{X_{\lambda_j}}(X_{\sigma_1} - X_0)$$

is a vector of martingales under $P^{\omega}_{z,0}$ with $E^{\omega}_{z,0}(M_k) = M_0 = z$. Let $\widetilde{M}_k$ denote the corresponding vector-valued martingale for $\widetilde{X}_k$. We have $E^{\omega}_{z,0}(\widetilde{M}_k) = \widetilde{M}_0 = 0$.

Let us observe that these martingales have nicely bounded moments. First by (3.9),

(5.4) $$\left| \sum_{j=0}^{\ell-1} \mathbb{1}\{X_{\lambda_j} \cdot \hat{u} = j\} E^{\omega}_{X_{\lambda_j}}(X_{\sigma_1} - X_0) \right| \leq C_3 \ell.$$

By another application of (3.9),

(5.5) $$E_{z,0}[(M_k - M_{k-1})^2 | \mathcal{H}_{k-1}] \leq C_4.$$

In particular, $M$ and $\widetilde{M}$ are $L^2$-martingales. We wish to apply optional stopping to the martingales $M$ and $\widetilde{M}$ and the stopping time $L$, justified by the next lemma. Given $\bar{p} \in (2, p)$, let us write

(5.6) $\hat{p} = (\bar{p} - 1) \wedge 2$, a number that satisfies $1 < \hat{p} \leq 2 < p$.

LEMMA 5.4. *Let $d \geq 1$ and consider a product probability measure $\mathbb{P}$ on environments with a forbidden direction $-\hat{u} \in \mathbb{Z}^d \setminus \{0\}$ as in (2.1). Assume nonnestling (N) in direction $\hat{u}$, and the moment hypothesis (M) with $p > 2$.*

*(a) There exists a constant $C_5$ such that $E_{z,0}(|\widetilde{X}_{\tilde{\lambda}_L}|^{\hat{p}}) \leq C_5$ for all choices of $z \in \mathbb{V}_d$. Same is true for $X_{\lambda_L} - z$.*

*(b) For $\mathbb{P}$-almost every $\omega$, $\{M_{\ell \wedge L} : \ell \geq 0\}$ and $\{\widetilde{M}_{\ell \wedge L} : \ell \geq 0\}$ are uniformly integrable martingales under $P^{\omega}_{z,0}$, for all choices of $z \in \mathbb{V}_d$.*

PROOF. Part (a). We do the proof for $\widetilde{X}_{\tilde{\lambda}_L}$. $\widetilde{M}_{L \wedge k}$ is also an $L^2$-martingale. By orthogonality of martingale increments, by $\mathcal{H}_{j-1}$-measurability of $\{L \geq j\} = \{L \leq j-1\}^c$, by (5.5), and by the integrability of $L$ (Lemma 5.3),

$$E_{z,0}(|\widetilde{M}_{L \wedge \ell}|^2) = \sum_{j=1}^{\ell} E_{z,0}(|\widetilde{M}_{L \wedge j} - \widetilde{M}_{L \wedge (j-1)}|^2)$$

$$= \sum_{j=1}^{\ell} E_{z,0}(|\widetilde{M}_j - \widetilde{M}_{j-1}|^2, L \geq j) \leq C_4 \sum_{j \geq 1} P_{z,0}(L \geq j) \leq C.$$

(Above $|\cdot|$ is the $\ell^2$-norm.) Then by Fatou's lemma $E_{z,0}(|\widetilde{M}_L|^2) \leq C$. Invoking (5.4) we finally get

$$E_{z,0}(|\widetilde{X}_{\widetilde{\lambda}_L}|^{\hat{p}}) \leq CE_{z,0}(|\widetilde{M}_L|^{\hat{p}}) + CE_{z,0}(L^{\hat{p}}) \leq C_5.$$

Part (b). The $\mathbb{P}$-full probability set of $\omega$'s is defined by the conditions $E^{\omega}_{z,0}(L^{\bar{p}-1}) < \infty$ and (3.9). This set is evidently of full $\mathbb{P}$-probability by Lemmas 3.1 and 5.3.

We prove the uniform integrability for $M_{\ell \wedge L}$, since the case of $\widetilde{M}_{\ell \wedge L}$ is the same. Due to (5.4), it suffices to check that $\{X_{\lambda(\ell \wedge L)}\}$ is uniformly integrable. By part (a) we only need to show the uniform integrability of $\{X_{\lambda_\ell} \mathbb{1}(L \geq \ell)\}$. For that, pick $q_1$ so that $1 < q_1 < \frac{1+\bar{p}}{2} \wedge (\bar{p}-1)$, let $q_2 = q_1/(q_1-1)$ be the conjugate exponent, and let $\nu = 1/q_2 = 1 - 1/q_1$. Then $q_1(1+\nu) = 2q_1 - 1 < \bar{p}$ and so (3.9) can be applied with exponent $q_1(1+\nu)$:

$$E^{\omega}_{z,0}[|X_{\lambda_\ell}|^{1+\nu} \mathbb{1}(L \geq \ell)]$$

$$\leq E^{\omega}_{z,0}\left[\ell^\nu \sum_{j=1}^{\ell} |X_{\lambda_j} - X_{\lambda_{j-1}}|^{1+\nu} \mathbb{1}(L \geq \ell)\right]$$

$$\leq \sum_{j=1}^{\infty} E^{\omega}_{z,0}[L^\nu \mathbb{1}(L \geq j)|X_{\lambda_j} - X_{\lambda_{j-1}}|^{1+\nu}]$$

$$\leq \sum_{j=1}^{\infty} (E^{\omega}_{z,0}[L])^{1/q_2} (E^{\omega}_{z,0}[\mathbb{1}(L \geq j)|X_{\lambda_j} - X_{\lambda_{j-1}}|^{(1+\nu)q_1}])^{1/q_1}$$

$$\leq C \sum_{j=1}^{\infty} (E^{\omega}_{z,0}[L])^{1/q_2} P^{\omega}_{z,0}(L \geq j)^{1/q_1}$$

$$\leq C \sum_{j=1}^{\infty} (E^{\omega}_{z,0}[L])^{1/q_2} (E^{\omega}_{z,0}[L^{\bar{p}-1}])^{1/q_1} j^{-(\bar{p}-1)/q_1} < \infty$$

where the convergence of the series comes from $(\bar{p}-1)/q_1 > 1$. In the second-to-last inequality we used the $\mathcal{H}_{\ell-1}$-measurability of the event $\{L \geq \ell\}$ and (3.9) with exponent $q_1(1+\nu)$. □

The conclusion from uniform integrability is that by optional stopping $E^{\omega}_{z,0}(M_L) = M_0$ and $E^{\omega}_{z,0}(\widetilde{M}_L) = \widetilde{M}_0$. With this we get

$$E^{\omega}_{z,0}(X_{\lambda_L} - \widetilde{X}_{\widetilde{\lambda}_L})$$

$$= E^{\omega}_{z,0}\left(M_L - \widetilde{M}_L + \sum_{j=0}^{L-1} \mathbb{1}\{X_{\lambda_j} \cdot \hat{u} = j\} E^{\omega}_{X_{\lambda_j}}(X_{\sigma_1} - X_0)\right)$$

$$-\sum_{j=0}^{L-1}\mathbb{1}\{\widetilde{X}_{\widetilde{\lambda}_j}\cdot\hat{u}=j\}E^\omega_{\widetilde{X}_{\widetilde{\lambda}_j}}(X_{\sigma_1}-X_0)\Bigg)$$

$$=z+E^\omega_{z,0}\Bigg(\sum_{j=0}^{L-1}\mathbb{1}\{X_{\lambda_j}\cdot\hat{u}=j\}E^\omega_{X_{\lambda_j}}(X_{\sigma_1}-X_0)$$

$$-\sum_{j=0}^{L-1}\mathbb{1}\{\widetilde{X}_{\widetilde{\lambda}_j}\cdot\hat{u}=j\}E^\omega_{\widetilde{X}_{\widetilde{\lambda}_j}}(X_{\sigma_1}-X_0)\Bigg).$$

Abbreviate

$$S=\sum_{j=0}^{L-1}\mathbb{1}\{X_{\lambda_j}\cdot\hat{u}=j\} \quad\text{and}\quad \widetilde{S}=\sum_{j=0}^{L-1}\mathbb{1}\{\widetilde{X}_{\widetilde{\lambda}_j}\cdot\hat{u}=j\}$$

for the numbers of levels that the walks visit before level $L$. Integrating out the environments then gives

$$E_{z,0}\Bigg[\sum_{j=0}^{L-1}\mathbb{1}\{X_{\lambda_j}\cdot\hat{u}=j\}E^\omega_{X_{\lambda_j}}(X_{\sigma_1}-X_0)\Bigg]$$

$$=\sum_{j=0}^{\infty}E_{z,0}[\mathbb{1}\{j<L\}\mathbb{1}\{X_{\lambda_j}\cdot\hat{u}=j\}E^\omega_{X_{\lambda_j}}(X_{\sigma_1}-X_0)]$$

$$=E_0(X_{\sigma_1})\sum_{j=0}^{\infty}E_{z,0}[\mathbb{1}\{j<L\}\mathbb{1}\{X_{\lambda_j}\cdot\hat{u}=j\}]$$

$$=E_0(X_{\sigma_1})E_{z,0}(S)$$

with a corresponding formula for the $\widetilde{X}$-walk. Substituting this back up leads to

$$E_{z,0}(X_{\lambda_L}-\widetilde{X}_{\widetilde{\lambda}_L})=z+E_0(X_{\sigma_1})E_{z,0}(S-\widetilde{S}).$$

Project this equation onto $\hat{u}$. Since $X_{\lambda_L}\cdot\hat{u}-\widetilde{X}_{\widetilde{\lambda}_L}\cdot\hat{u}=0$ by the definition of $L$ and $z\cdot\hat{u}=0$ while $X_{\sigma_1}\cdot\hat{u}\geq 1$, we conclude that $E_{z,0}(S-\widetilde{S})=0$. Substituting this back up gives this conclusion:

(5.7) $$E_{z,0}(X_{\lambda_L}-\widetilde{X}_{\widetilde{\lambda}_L})=z,$$

which is a mean-zero increment property.

Recall the definition (5.1) of the $\mathbb{V}_d$-valued Markov chain $Z_n$ that tracks the difference of the walks $X$ and $\widetilde{X}$ on successive new common levels. The transition probability of $Z_n$ is given for $x,y\in\mathbb{V}_d$ by

$$q(x,y)=P_{x,0}[X_{\lambda_L}-\widetilde{X}_{\widetilde{\lambda}_L}=y].$$

To paraphrase the formula, in order to find the next state $y$ from the present state $x$, put the $X$-walk at $x$, put the $\widetilde{X}$-walk at the origin, let the walks run until they have reached a new common level $L$ above 0 and let $y$ be the difference of the entry points at level $L$.

We are now all set for controlling the chain $(Z_k)$. Recall that $\hat{p} = (\bar{p}-1) \wedge 2$ and $\bar{p} \in (2, p)$.

LEMMA 5.5. *Let $d \geq 1$ and consider a product probability measure $\mathbb{P}$ on environments with a forbidden direction $-\hat{u} \in \mathbb{Z}^d \setminus \{0\}$ as in (2.1). Assume nonnestling* (N) *in direction $\hat{u}$ and the moment hypothesis* (M) *with $p > 2$. Then the transition $q(x, y)$ has these properties for all $x \in \mathbb{V}_d$:*

$$\sum_{m \in \mathbb{V}_d} m q(x, x+m) = 0 \tag{5.8}$$

*and there exists a constant $C_6 < \infty$ such that*

$$\sum_{m \in \mathbb{V}_d} |m|^{\hat{p}} q(x, x+m) \leq C_6. \tag{5.9}$$

*In addition to the assumptions above, assume $d \geq 2$ and ellipticity* (E). *Then there exists a constant $\varepsilon > 0$ such that*

$$q(x, x) \leq 1 - \varepsilon \qquad \text{for all } x \in \mathbb{V}_d. \tag{5.10}$$

PROOF. Property (5.8) follows from (5.7), and property (5.9) from Lemma 5.4(a).

We prove property (5.10) carefully, for even though the argument is elementary, it is here that the proof needs the ellipticity hypotheses. By assumption (2.3) in the ellipticity hypothesis (E) we can fix two nonzero vectors $z \neq y$ such that $\mathbb{E}(\pi_{0z} \pi_{0y}) > 0$. Pick their names so that $z \cdot \hat{u} \geq y \cdot \hat{u}$. If $y \cdot \hat{u} = z \cdot \hat{u} = 0$, then by nonnestling (N) there exists a vector $u$ with $u \cdot \hat{u} > 0$ and $\mathbb{E}(\pi_{0z} \pi_{0y} \pi_{0u}) > 0$. Thus by replacing $z$ with $u$ if necessary we can assume $z \cdot \hat{u} > 0$. Recall that we are assuming $\hat{u}$ is an integer vector, so all the dot products are also integers.

Let $x \in \mathbb{V}_d$. We distinguish three cases.

*Case 1.* $y \cdot \hat{u} = 0$. Then

$$1 - q(x,x) \geq q(x, x-y) \geq P_{x,0}(X_1 = x+z, \widetilde{X}_1 = y, \widetilde{X}_2 = y+z)$$

$$= \mathbb{E} \pi_{x, x+z} \pi_{0, y} \pi_{y, y+z} = \begin{cases} (\mathbb{E} \pi_{0z})^2 \mathbb{E} \pi_{0y}, & \text{if } x \notin \{0, y\}, \\ \mathbb{E}[\pi_{0z} \pi_{0y}] \mathbb{E} \pi_{0z}, & \text{if } x = 0, \\ \mathbb{E}[\pi_{0z}^2] \mathbb{E} \pi_{0y}, & \text{if } x = y. \end{cases}$$

*Case 2.* $y \cdot \hat{u} > 0$ and $y \notin \mathbb{R}z$. Let $n, m \geq 1$ be such that $nz \cdot \hat{u} = my \cdot \hat{u}$ is the least common multiple of $y \cdot \hat{u}$ and $z \cdot \hat{u}$. We have

$$1 - q(x, x)$$

$$\geq q(x, x + nz - my)$$
$$\geq P_{x,0}(X_i - X_{i-1} = z, \widetilde{X}_j - \widetilde{X}_{j-1} = y, \text{ for } i = 1, \ldots, n \text{ and } j = 1, \ldots, m)$$
$$= \begin{cases} (\mathbb{E}\pi_{0z})^n (\mathbb{E}\pi_{0y})^m, & \text{if } x \neq 0, \\ \mathbb{E}[\pi_{0z}\pi_{0y}](\mathbb{E}\pi_{0z})^{n-1}(\mathbb{E}\pi_{0y})^{m-1}, & \text{if } x = 0. \end{cases}$$

*Case* 3. $y \cdot \hat{u} > 0$ and $y \in \mathbb{R}z$. Together with the earlier assumption $y \cdot \hat{u} \leq z \cdot \hat{u}$ these imply $y \cdot \hat{u} < z \cdot \hat{u}$. The ellipticity hypothesis (E) implies the existence of a vector $w \notin \mathbb{R}z$ such that $\mathbb{E}(\pi_{0w}) > 0$. Consider the positive integer solutions $(\ell, m, n)$ of the equation

$$\ell(z \cdot \hat{u}) = m(y \cdot \hat{u}) + n(w \cdot \hat{u}).$$

Such solutions exist. For if $w \cdot \hat{u} = 0$, then $\ell = y \cdot \hat{u}$, $m = z \cdot \hat{u}$ together with any $n > 0$ works. If $w \cdot \hat{u} > 0$, then one solution is $\ell = w \cdot \hat{u}$, $m = w \cdot \hat{u}$, $n = z \cdot \hat{u} - y \cdot \hat{u}$. Fix a solution where $\ell$ is minimal. Define a path $(\tilde{x}_j)_{j=1}^{m+n}$ such that $\tilde{x}_1 = y$, $\tilde{x}_2 = y + w$, and after that each step is either $y$ or $w$ but so that $\tilde{x}_{m+n} = my + nw$. Define another path $(x_k = kz)_{k=1}^{\ell}$. Paths $(\tilde{x}_j)_{j=1}^{m+n}$ and $(x_k)_{k=1}^{\ell}$ do not have a common level until at $\tilde{x}_{m+n} \cdot \hat{u} = x_\ell \cdot \hat{u}$. To see this, note two points:

(i) $\tilde{x}_1 \cdot \hat{u} = x_k \cdot \hat{u}$ is impossible due to the assumption $y \cdot \hat{u} < z \cdot \hat{u}$.
(ii) An equality $\tilde{x}_j \cdot \hat{u} = x_k \cdot \hat{u}$ with $2 \leq j \leq m + n$ and $1 \leq k < \ell$ would produce a solution $(\ell, m, n)$ with smaller $\ell$.

Note also that, since $z, y, w \neq 0$ and by the linear independence of $w$ and $z$,

$$x + x_\ell - \tilde{x}_{m+n} = x + \ell z - my - nw \neq x \qquad \text{for any } x.$$

So
$$1 - q(x,x) \geq q(x, x + x_\ell - \tilde{x}_{m+n})$$
$$\geq P_{x,0}(X_k = x + x_k \text{ for } 1 \leq k \leq \ell \text{ and } \widetilde{X}_j = \tilde{x}_j \text{ for } 1 \leq j \leq m+n)$$
$$= \begin{cases} (\mathbb{E}\pi_{0z})^\ell (\mathbb{E}\pi_{0y})^m (\mathbb{E}\pi_{0w})^n, & \text{if } x \neq 0, \\ \mathbb{E}[\pi_{0z}\pi_{0y}](\mathbb{E}\pi_{0z})^{\ell-1}(\mathbb{E}\pi_{0y})^{m-1}(\mathbb{E}\pi_{0w})^n, & \text{if } x = 0. \end{cases}$$

The three cases give finitely many positive lower bounds on $1 - q(x,x)$ that are independent of $x$. Let $\varepsilon$ be their minimum. $\square$

We can now finish the proof of Proposition 5.1. Write

$$G_n(x,y) = \sum_{k=0}^{n} q^k(x,y) = E_x\left[\sum_{k=0}^{n} \mathbb{1}\{Z_k = y\}\right],$$

where $q^k(x,y)$ is the $k$-step transition probability from $x$ to $y$, and now $E_x$ is the expectation on the path space of $\{Z_k\}$ when the initial state is $Z_0 = x$.

Continue from (5.2) and apply Lemma A.4 from the Appendix and the summability of $h$:

$$E_{0,0}(|X_{[0,n-1]} \cap \widetilde{X}_{[0,n-1]}|) \leq \sum_{j=0}^{n-1} E_0[h(Z_j)] = \sum_{x \in \mathbb{V}_d} h(x) G_{n-1}(0,x)$$
$$\leq C_6 n^{1/\hat{p}} \sum_{x \in \mathbb{V}_d} h(x) \leq C n^{1/\hat{p}}.$$

Recalling that $\hat{p} = (\bar{p} - 1) \wedge 2$ and $\bar{p} \in (2, p)$, this completes the proof of Proposition 5.1, and thereby the proof of Theorem 2.1. □

## APPENDIX

**A.1. A linear algebra lemma.** Let us say that a vector is rational if it has rational coordinates, and analogously a vector is integral if it has integer coordinates. The lemma below implies that for a finitely supported step distribution the requirement of a rational vector in the forbidden direction assumption (2.1) and in the nonnestling hypothesis (N) is no more restrictive than requiring a general real vector. This justifies the derivation of Corollary 2.2 from Theorem 2.1.

LEMMA A.1. *Let $A$ be a finite subset of $\mathbb{Z}^d$. Suppose there exists a vector $\hat{v} \in \mathbb{R}^d$ such that $\hat{v} \cdot x \geq 0$ for all $x \in A$. Then there exists a vector $\hat{u}$ with integer coordinates such that, for all $x \in A$,*

$$\hat{u} \cdot x > 0 \quad \text{if and only if} \quad \hat{v} \cdot x > 0 \quad \text{and}$$
$$\hat{u} \cdot x = 0 \quad \text{if and only if} \quad \hat{v} \cdot x = 0.$$

The proof is mainly done in the next lemma. Let us recall this notion of vector product: if $h_1, \ldots, h_{d-1}$ are vectors in $\mathbb{R}^d$, let $z = F(h_1, \ldots, h_{d-1})$ be the vector defined by the equations

$$\det[h_1, \ldots, h_{d-1}, x] = x \cdot z, \qquad \text{for all } x \in \mathbb{R}^d.$$

Here $[h_1, \ldots, h_{d-1}, x]$ denotes a matrix in terms of its column decomposition. Explicitly, $z = [z(1), \ldots, z(d)]^t$ with coordinates

$$z(i) = (-1)^{i+d} \det[h_1, \ldots, h_{d-1}]\{i\},$$

where $[h_1, \ldots, h_{d-1}]\{i\}$ is the $(d-1) \times (d-1)$ matrix obtained from $[h_1, \ldots, h_{d-1}]$ by removing row $i$. Consequences of the definition are that $z \cdot h_i = 0$ for each $h_i$, and $z \neq 0$ if, and only if, $h_1, \ldots, h_{d-1}$ are linearly independent. The explicit formula shows that if all $h_i$ are integer (resp. rational) vectors, then so is $z$.

LEMMA A.2. *Let $v_1, \ldots, v_n$ be linearly independent vectors in $\mathbb{R}^d$ that lie in the orthogonal complement $\{\hat{v}\}^\perp$ of some vector $\hat{v} \in \mathbb{R}^d$. Suppose $v_1, \ldots, v_n$ all have integer coordinates. Then for each $\varepsilon > 0$ there exists a vector $w$ with rational coordinates such that $|w - \hat{v}| \leq \varepsilon$ and $v_1, \ldots, v_n \in \{w\}^\perp$.*

PROOF. If $n = d - 1$, then $z = F(v_1, \ldots, v_{d-1})$ is a vector with integer coordinates and the property $\operatorname{span}\{v_1, \ldots, v_{d-1}\} = \{z\}^\perp$. The spans of $z$ and $\hat{v}$ must then coincide, so in particular we can take rational multiples of $z$ arbitrarily close to $\hat{v}$.

Assume now $n < d - 1$. Find vectors $\xi_{n+1}, \ldots, \xi_{d-1}$ so that

$$v_1, \ldots, v_n, \xi_{n+1}, \ldots, \xi_{d-1}$$

is a basis for $\{\hat{v}\}^\perp$. Next find rational vectors $\eta_{n+1}^m, \ldots, \eta_{d-1}^m$ such that for each $n+1 \leq k \leq d-1$, $\eta_k^m \to \xi_k$ as $m \to \infty$, and so that

$$v_1, \ldots, v_n, \eta_{n+1}^m, \ldots, \eta_{d-1}^m$$

are linearly independent for each $m$.

This can be achieved by the following argument. Suppose that for a particular $m \geq 1$ and $n \leq k < d-1$, vectors $\eta_{n+1}^m, \ldots, \eta_k^m$ have been chosen so that $|\eta_j^m - \xi_j| \leq 1/m$ for $n+1 \leq j \leq k$, and the system $v_1, \ldots, v_n, \eta_{n+1}^m, \ldots, \eta_k^m$ is linearly independent. The case $k = n$ corresponds to the case where none of these vectors has been chosen yet, for the given $m$. The subspace

$$U = \operatorname{span}\{v_1, \ldots, v_n, \eta_{n+1}^m, \ldots, \eta_k^m\}$$

has dimension $k < d - 1$ and is a closed subset with empty interior in $\mathbb{R}^d$. Consequently the set $B_{1/m}(\xi_{k+1}) \setminus U$ is nonempty and open, and we can choose any rational vector $\eta_{k+1}^m$ from this set.

Once the rational vectors $\eta_{n+1}^m, \ldots, \eta_{d-1}^m$ have been defined, let

$$\zeta_m = F(v_1, \ldots, v_n, \eta_{n+1}^m, \ldots, \eta_{d-1}^m).$$

This $\zeta_m$ is a rational vector. Next let $s_m$ be the real number defined by

$$|s_m \zeta_m - \hat{v}| = \inf\{|t\zeta_m - \hat{v}| : t \in \mathbb{R}\}$$

and then let $q_m$ be a rational such that $|s_m - q_m| < 1/m$. Finally, let $w_m = q_m \zeta_m$. Clearly, $\{w_m\}^\perp$ contains $v_1, \ldots, v_n$. We claim that $w_m \to \hat{v}$ as $m \to \infty$.

The product $F$ is continuous in its arguments, so

$$\zeta_m \to \zeta = F(v_1, \ldots, v_n, \xi_{n+1}, \ldots, \xi_{d-1}).$$

Since $\zeta$ and $\hat{v}$ both span the orthogonal complement of $\{v_1, \ldots, v_n, \xi_{n+1}, \ldots, \xi_{d-1}\}$, there is a real $s$ such that $\hat{v} = s\zeta$. Consequently $s\zeta_m \to s\zeta = \hat{v}$. Now

$$|w_m - \hat{v}| = |q_m \zeta_m - \hat{v}| \leq |q_m - s_m| \cdot |\zeta_m| + |s_m \zeta_m - \hat{v}|.$$

The first term after the above inequality vanishes as $m \to \infty$ by the choice of $q_m$ and because $|\zeta_m| \to |\zeta|$. By the definition of $s_m$

$$|s_m \zeta_m - \hat{v}| \leq |s\zeta_m - \hat{v}| \to 0$$

as observed earlier. This completes the proof of the lemma. $\square$

PROOF OF LEMMA A.1. Let

$$M = \max\{|x| : x \in A\} < \infty \quad \text{and} \quad \delta = \min\{\hat{v} \cdot x : x \in A, \hat{v} \cdot x > 0\} > 0.$$

Let $v_1, \ldots, v_n$ be a maximal linearly independent set from $A \cap \{\hat{v}\}^\perp$. If this set is not empty, then pick a rational vector $w$ from Lemma A.2 with $\varepsilon = \delta/(2M)$. Otherwise, just pick any rational vector $w \in B_\varepsilon(\hat{v})$. Then for $x \in A$ we have on the one hand

$$\hat{v} \cdot x = 0 \Longrightarrow x \in \text{span}\{v_1, \ldots, v_n\} \Longrightarrow w \cdot x = 0,$$

and on the other hand

$$\hat{v} \cdot x > 0 \Longrightarrow \hat{v} \cdot x \geq \delta \Longrightarrow w \cdot x \geq \hat{v} \cdot x - |(w - \hat{v}) \cdot x| \geq \delta - M|w - \hat{v}| \geq \delta/2.$$

Now let $\hat{u}$ be a large enough positive integer multiple of $w$. $\square$

**A.2. A renewal process bound.** Write $\mathbb{N}^* = \{1, 2, 3, \ldots\}$ and $\mathbb{N} = \{0, 1, 2, \ldots\}$. The setting for the next technical lemma is the following. Let $\{Y_i : i \in \mathbb{N}^*\}$ be a sequence of i.i.d. positive integer-valued random variables, and $\{\widetilde{Y}_j : j \in \mathbb{N}^*\}$ an independent copy of this sequence. $Y$ denotes a random variable with the same distribution as $Y_1$. The corresponding renewal processes are defined by

$$S_0 = \widetilde{S}_0 = 0, S_n = Y_1 + \cdots + Y_n \quad \text{and} \quad \widetilde{S}_n = \widetilde{Y}_1 + \cdots + \widetilde{Y}_n \qquad \text{for } n \geq 1.$$

Let $h$ be the largest positive integer such that the common distribution of $Y_i$ and $\widetilde{Y}_j$ is supported on $h\mathbb{N}^*$. For $i, j \in h\mathbb{N}$ define

$$L_{i,j} = \inf\{\ell \geq 1 : \text{there exist } m, n \geq 0 \text{ such that } i + S_m = \ell = j + \widetilde{S}_n\}.$$

The restriction $\ell \geq 1$ in the definition has the consequence that $L_{i,i} = i$ for $i > 0$ but $L_{0,0}$ is nontrivial. The next lemma is proved in the Appendix of [11].

LEMMA A.3. *Let $1 \leq r < \infty$ be a real number, and assume $E(Y^{r+1}) < \infty$. Then there exists a finite constant $C$ such that for all $i, j \in h\mathbb{N}$,*

$$E(L_{i,j}^r) \leq C(1 + i^r + j^r).$$

### A.3. A Green function estimate.

LEMMA A.4. *Let $\mathbb{V}$ be a subset of some $\mathbb{Z}^d$, $d \geq 1$. Consider a Markov chain $Z_n$ on $\mathbb{V}$ whose transition $q(x,y)$ has properties* (5.8)–(5.10) *with $1 < \hat{p} \leq 2$. Then there exists a constant $C_7 < \infty$ such that*

$$G_n(x,y) = \sum_{k=0}^{n} q^k(x,y) = E_x\left(\sum_{k=0}^{n} \mathbb{1}\{Z_k = y\}\right) \leq C_7 n^{1/\hat{p}}$$

*for all $n \geq 1$ and all $x, y \in \mathbb{V}$.*

PROOF. First we use the familiar argument to reduce the proof to the diagonal case. For $k \geq 1$, let

$$f^k(x,y) = P_x(Z_1 \neq y, \ldots, Z_{k-1} \neq y, Z_k = y)$$

be the probability that after time 0 the first visit from $x$ to $y$ occurs at time $k$. Note that $\sum_k f^k(x,y) \leq 1$. Then for $x \neq y$

$$G_n(x,y) = \sum_{k=1}^{n} q^k(x,y) = \sum_{k=1}^{n}\sum_{j=1}^{k} f^j(x,y) q^{k-j}(y,y)$$

$$= \sum_{j=1}^{n} f^j(x,y) \sum_{k=j}^{n} q^{k-j}(y,y)$$

$$\leq \sum_{j=1}^{n} f^j(x,y) \sum_{k=0}^{n} q^k(y,y) \leq \sum_{k=0}^{n} q^k(y,y) = G_n(y,y).$$

We can now take $x = y$ and it remains to show

$$E_x\left(\sum_{k=0}^{n} \mathbb{1}\{Z_k = x\}\right) \leq C_7 n^{1/\hat{p}}.$$

Keep $x$ fixed now, and consider the Markov chain $Z_n$ under the measure $P_x$ on its path space. By properties (5.8) and (5.9), $Z_n$ is an $L^{\hat{p}}$-martingale relative to its own filtration $\{\mathcal{F}_n^Z\}$, with initial point $Z_0 = x$. Furthermore, (5.9) implies a uniform bound on conditional $\hat{p}$th moments of increments:

(A.1) $$E_x(|Z_k - Z_{k-1}|^{\hat{p}}|\mathcal{F}_{k-1}^Z) \leq C_6.$$

Let $0 = \tau_0 < \tau_1 < \tau_2 < \cdots$ be the successive times of arrivals to $x$ after leaving $x$, in other words

$$\tau_{j+1} = \inf\{n > \tau_j : Z_n = x \text{ and } Z_k \neq x \text{ for some } k : \tau_j < k < n\}.$$

Let $T_j$ ($j \geq 0$) be the durations of the sojourns at $x$; in other words

$$Z_n = x \quad \text{if, and only if} \quad \tau_j \leq n < \tau_j + T_j \quad \text{for some } j \geq 0.$$

Given that an arrival has happened, the sojourns are independent of the past and have geometric distributions, so on the event $\{\tau_j < \infty\}$,

$$E_x(T_j | \mathcal{F}^Z_{\tau_j}) = \frac{1}{1 - q(x,x)}.$$

Let $J_n = \max\{j \geq 0 : \tau_j \leq n\}$ mark the last sojourn at $x$ that started by time $n$. With these notations

$$E_x\left(\sum_{k=0}^n \mathbb{1}\{Z_k = x\}\right) \leq E_x\left(\sum_{j=0}^{J_n} T_j\right)$$

(A.2)
$$= \sum_{j=0}^\infty E_x(\mathbb{1}\{\tau_j \leq n\} T_j)$$

$$= \frac{1}{1 - q(x,x)} E_x(1 + J_n).$$

To bound the number $J_n$ of arrivals to $x$ from somewhere else we use the upcrossing lemma from martingale theory. Write $Z_n = (\xi^1_n, \ldots, \xi^d_n)$ in terms of the (standard) coordinates, and similarly $x = (t^1, \ldots, t^d)$. Let $U^i_n$ count the number of upcrossings of the martingale $\xi^i$ across the interval $[t^i - 1, t^i]$ up to time $n$. Similarly $V^i_n$ counts the number of downcrossings across the interval $[t^i, t^i + 1]$ made by the martingale $\xi^i$ up to time $n$. Quite obviously

$$J_n \leq \sum_{i=1}^d (U^i_n + V^i_n)$$

since each arrival to $x$ means that some coordinate arrived at $t^i$ from either above or below. By the upcrossing inequality ([6], (2.9) in Chapter 4)

$$E_x(U^i_n) \leq E_x[(\xi^i_n - (t^i - 1))^+] - E_x[(\xi^i_0 - (t^i - 1))^+]$$
$$\leq E_x[|\xi^i_n - t^i|] + 1 - 1 = E_x[|\xi^i_n - t^i|].$$

Similarly $E_x(V^i_n) \leq E_x[|\xi^i_n - t^i|]$, by applying the upcrossing inequality to $-\xi^i$ and the interval $[-t_i - 1, -t_i]$. Now follows

$$E_x[J_n] \leq \sum_{i=1}^d (E_x[U^i_n] + E_x[V^i_n])$$

$$\leq 2 \sum_{i=1}^d E_x[|\xi^i_n - t^i|] = 2 \sum_{i=1}^d E_x[|\xi^i_n - \xi^i_0|].$$

In the next stage we apply the increment bound (A.1). Since $1 < \hat{p} \leq 2$ we can apply the Burkholder–Davis–Gundy inequality ([4], Theorem 3.2) to

derive

$$2\sum_{i=1}^{d} E_x[|\xi_n^i - \xi_0^i|] \leq 2\sum_{i=1}^{d} \{E_x[|\xi_n^i - \xi_0^i|^{\hat{p}}]\}^{1/\hat{p}}$$
$$\leq 2C \sum_{i=1}^{d} \left\{ E_x \left[ \left( \sum_{k=1}^{n} (\xi_k^i - \xi_{k-1}^i)^2 \right)^{\hat{p}/2} \right] \right\}^{1/\hat{p}}$$
$$\leq 2C \sum_{i=1}^{d} \left\{ E_x \sum_{k=1}^{n} |\xi_k^i - \xi_{k-1}^i|^{\hat{p}} \right\}^{1/\hat{p}}$$
$$\leq C n^{1/\hat{p}}.$$

The next-to-last inequality came from noticing that $\hat{p}/2 \in (0,1]$ and hence, for any nonnegative summands,

$$(x_1 + \cdots + x_n)^{\hat{p}/2} \leq x_1^{\hat{p}/2} + \cdots + x_n^{\hat{p}/2}.$$

Substituting the bounds back up to line (A.2) and applying property (5.10) gives

$$G_n(x,x) = E_x\left[\sum_{k=0}^{n} \mathbb{1}\{Z_k = x\}\right] \leq \frac{1 + Cn^{1/\hat{p}}}{1 - q(x,x)} \leq \frac{1 + Cn^{1/\hat{p}}}{\varepsilon} \leq C_7 n^{1/\hat{p}}$$

for a new constant $C_7$. The proof of the lemma is complete. $\square$

**Acknowledgment.** The authors thank an anonymous referee for a careful reading of the paper.

DEPARTMENT OF MATHEMATICS
UNIVERSITY OF UTAH
155 SOUTH 1400 EAST
SALT LAKE CITY, UTAH 84112
USA
E-MAIL: firas@math.utah.edu
URL: www.math.utah.edu/~firas

DEPARTMENT OF MATHEMATICS
UNIVERSITY OF WISCONSIN–MADISON
VAN VLECK HALL
MADISON, WISCONSIN 53706
USA
E-MAIL: seppalai@math.wisc.edu
URL: www.math.wisc.edu/~seppalai